\documentclass{article}
\usepackage{spconf,amsmath,graphicx}
\usepackage{mathtools}
\usepackage[ruled,vlined,linesnumbered]{algorithm2e}
\usepackage{dsfont}
\usepackage{tikz}
\DeclareMathOperator*{\argmin}{arg\,min}

\title{A temporal multiscale approach for MR Fingerprinting}
\name{Samuel Cortinhas$^{\star}$, Mohammad Golbabaee$^{\dagger}$ and Matthias J. Ehrhardt$^{\times}$}
\address{${}^{\star}$Department of Mathematical Sciences, ${}^{\dagger}$Department of Computer Science, \\ ${}^{\times}$Institute for Mathematical Innovation, University of Bath, UK}
\begin{document}
\ninept
\maketitle
\begin{abstract}
Quantitative MRI (qMRI) is becoming increasingly important for research and clinical applications, however, state-of-the-art reconstruction methods for qMRI are computationally prohibitive. We propose a temporal multiscale approach to reduce computation times in qMRI. Instead of computing exact gradients of the qMRI likelihood, we propose a novel approximation relying on the temporal smoothness of the data. These gradients are then used in a coarse-to-fine (C2F) approach, for example using coordinate descent. The C2F approach was also found to improve the accuracy of solutions, compared to similar methods where no multiscaling was used.
\end{abstract}
\begin{keywords}
quantitative MRI, MR fingerprinting, nonlinear inverse problems, coordinate descent.
\end{keywords}
\section{Introduction}
\label{sec:intro}

Quantitative MRI improves on standard (qualitative) MRI by providing a map of tissue parameters, e.g. proton density ($\rho$), longitudinal ($T_1$) and transverse ($T_2$) relaxation times, and off-resonance frequency ($\omega$). These parameters are related to the physics of nuclear magnetic resonance via the Bloch equations \cite{bloch1946nuclear}. Quantifying this information has the potential benefits of being more sensitive to early indicators of disease, among others \cite{ma2013magnetic}. 

Magnetic Resonance Fingerprinting (MRF) \cite{ma2013magnetic}, is a state-of-the-art technology to make multi-parametric qMRI acquisitions fast. The proposed numerical approaches for MRF \cite{davies2014compressed, mcgivney2014svd, asslander2018low} are mostly based on template-matching signals to a presimulated dictionary of possible solutions discretised over a domain. However, limitations of this approach include memory problems and low accuracy.

Dictionary-free reconstruction using a Bloch response model \cite{dong2019quantitative, sbrizzi2018fast} requires solving a nonlinear fitting problem. For this, computing gradients requires evaluations of the physical forward mapping (i.e. solving the Bloch equations) and their derivatives across a potentially very long temporal sequence, which can severely slow down computations. In this paper, we introduce a multiscale approach across the temporal domain to reduce computation times in qMRI. Furthermore, we present a C2F approach using a variant of coordinate descent.

\subsection{MRF Signal model}
\label{ssec:model}

We consider an Inversion Recovery balance Steady-State Free Precession (IR-bSSFP) excitation pulse sequence \cite{scheffler1999pictorial} of length $L$, with constant repetition time $T_R$ (e.g. $10$ ms). The magnetisation $m_{l}$ after the $l$-th excitation pulse, measured at the middle of the $T_R$ time interval, satisfies the recursive formula
\begin{equation}
m_{l}(u)=A_{l} m_{l-1}(u)+b
\end{equation}
where $u=(T_1,T_2,\omega)$, $A_{l}=E V_{\omega} R_{\alpha_{l}} V_{\omega}^{T}$ and $b=(1-e_2) m_{e}$ \cite{dong2019quantitative, scheffler1999pictorial}. In particular, $m_{e}=(0,0,1)^{T}$, $m_{0}=-m_e$, $e_1= \exp(-T_R/T_2)$, $e_2= \exp(-T_R/T_1)$, $E=\operatorname{diag}(e_1,e_1,e_2)$. Finally, $(\alpha_{l})_{l=1}^{L}$ denotes the flip angle sequence, $\phi=2 \pi \omega T_R$ and $R_{\alpha_l}$, $V_{\omega}$ are 3D rotation matrices (see \cite{dong2019quantitative}).

We follow the modelling of Dong et al. \cite{dong2019quantitative}, which consists of the nonlinear equation $Q_{l}(\rho,u)=y_{l}$, where $y_l$ denotes the measured $k$-space data and the nonlinear forward operator is defined by
\begin{equation}
Q_{l}(\rho,u)=P_{l} \mathcal{F} \left(\rho T_{x,y} m_{l}\left(u\right)\right),
\end{equation}
where $m_l$ is defined via (1), $T_{x,y}$ is the transverse projection onto the complex plane, i.e. $T_{x,y}m_{l}=[m_{l}]_{x}+i [m_{l}]_{y}$, $\mathcal{F}$ is the Fourier transform and $P_l$ is a temporally-varying $k$-space subsampling operator. 

\section{Proposed multiscale approach}
\label{sec:multiscaling}

\subsection{Multiscale Bloch mapping}

Instead of computing the full Bloch response $(m_{l})_{l=1}^{L}$, we seek an efficient method to approximate a subset of this sequence, i.e. $(m_{l})_{l \in S}$, with $S \subseteq \{1,2,..., L\}$. For convenience, we introduce the notation
\begin{equation*}
S_N(\delta)=\{\delta,\delta+ N, \delta + 2 N,..., \delta+(\lfloor L/N \rfloor -1) N\},
\end{equation*}
to denote a uniform grid with increment size $N$ and offset $\delta \in \{1,2,...,N\}$.

\begin{figure}

\begin{minipage}[b]{1.0\linewidth}
  \centering
  \centerline{\includegraphics[width=8.5cm,trim={3cm 0.2cm 3.5cm 0.3cm},clip]{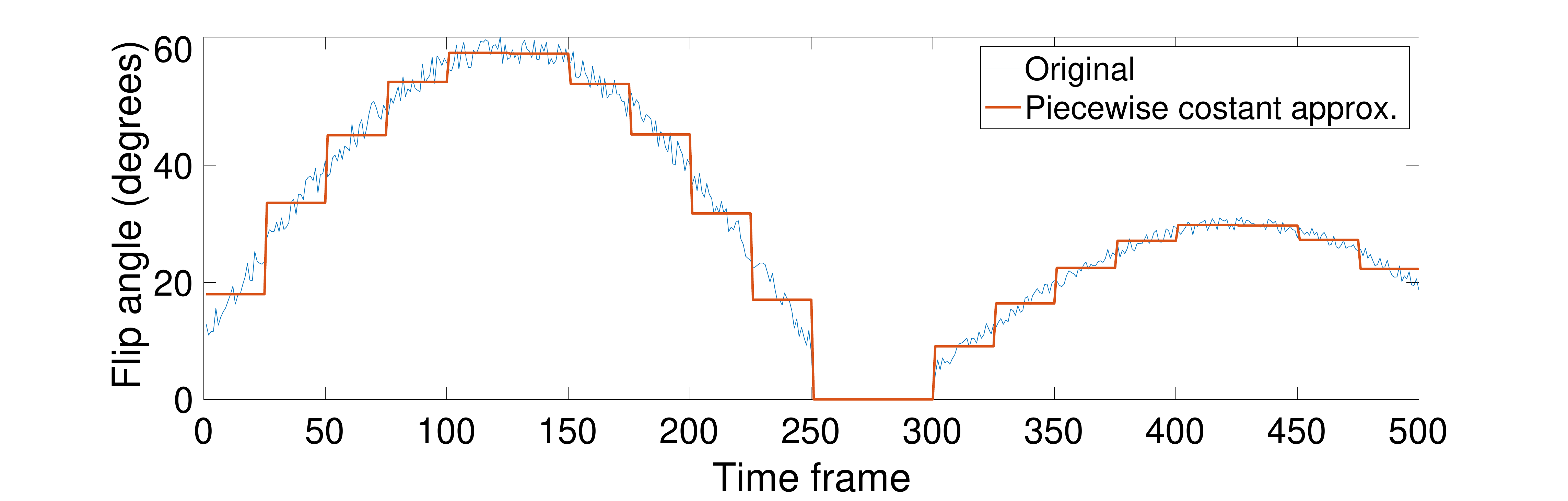}}
%  \vspace{2.0cm}
  \centerline{(a) Flip angle sequences}\medskip
\end{minipage}
\begin{minipage}[b]{1.0\linewidth}
  \centering
  \centerline{\includegraphics[width=8.5cm,trim={2.5cm 0.2cm 3.5cm 0.3cm},clip]{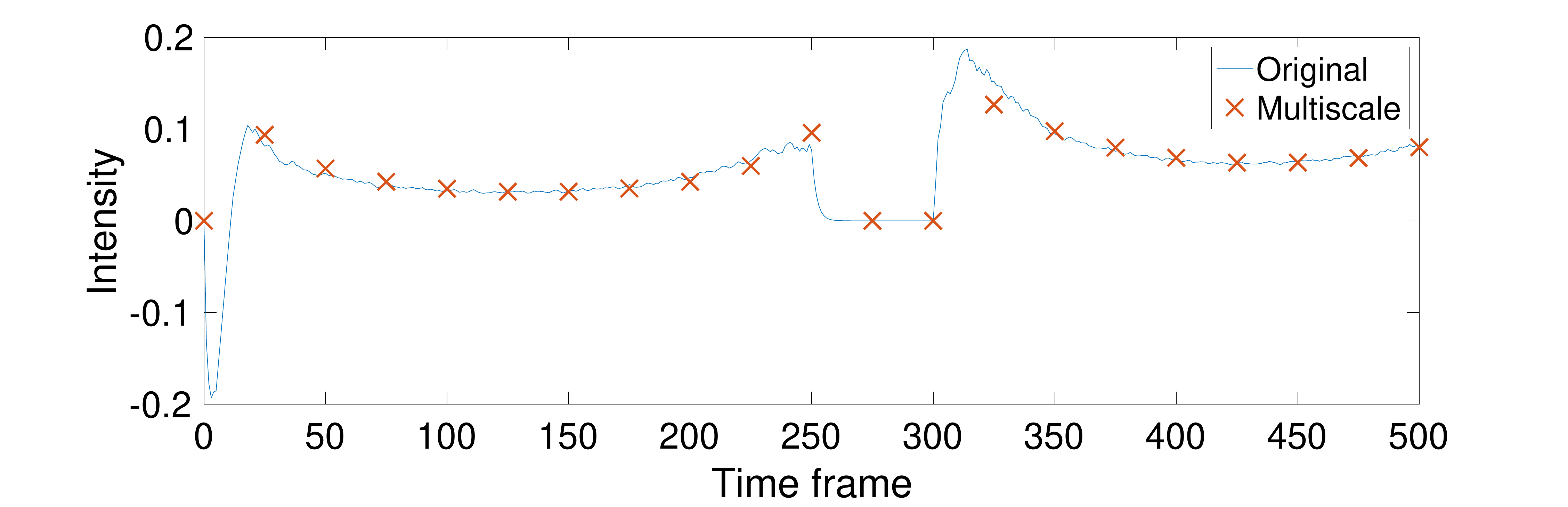}}
%  \vspace{1.5cm}
  \centerline{(b) Bloch responses ($y$-component)}
\end{minipage}

\caption{(a) Comparison of a typical flip angle sequence used in qMRI \cite{ma2013magnetic} with a piecewise constant approximation over the uniform grid $S_{25}(25)$ and (b) the resulting Bloch responses with tissue parameters $u=(811,77,0)$ (ms).}
\label{fig:flip}
\end{figure}

Note that the dependency of $A_{l}$ on $l$ (see (1)) only appears as a result of the flip angle sequence $(\alpha_{l})_{l=1}^{L}$. With this in mind, we make the simplifying assumption that $\alpha_{i}=\overline{\alpha}$, $i \in \{l+1,...,l+N\}$, with $\overline{\alpha}$ chosen to be the mean of the $\alpha_{i}$ over this interval (see fig. 1 for an example). The iteration matrices $A_{l}$ then lose their dependency on $l$ over each interval, i.e. $A_{i}=A$, $i \in \{l+1,...,l+N\}$, and (1) can be simplified as
\begin{equation}
m_{l+N}=A^{N} m_{l}+ \sum_{k=0}^{N-1} A^{k} b.
\end{equation}

For an efficient implementation of (3) we rely on the spectral decomposition of $A$. Its eigenvalues $\lambda_1$, $\lambda_2$, $\lambda_3$ and corresponding eigenvectors $v_1$, $v_2$, $v_3$ are easily found by observing the similarity transform $A=V_{\omega} \left(E R_{\overline{\alpha}}\right) V_{\omega}^{-1}$. We then define $D=\operatorname{diag}(\lambda_{1}, \lambda_{2}, \lambda_{3})$ and $Z=[v_{1},v_{2},v_{3}]$ to write $A=Z D Z^{-1}$. In particular, for $(G_i)_{i=1}^3$ denoting the $3$$\times$$3$ elementary matrices with a $1$ in the $i$-th diagonal slot and $0$'s elsewhere, we have
\begin{equation}
        A^{k}=\sum_{i=1}^{3} \lambda_{i}^{k} U_{i}
\end{equation}
where $U_{i}=Z G_{i} Z^{-1}$ for $i = 1,2,3$ and $k \in \mathds{N}$. Putting (4) into~(3) leads to the formula 
\begin{equation*}
    m_{l+N}=\left(\sum_{i=1}^{3} \lambda_{i}^{N} U_{i} \right) m_{l} + \left(\sum_{i=1}^{3} \frac{1-\lambda_{i}^{N}}{1-\lambda_{i}} U_{i} \right) b.
\end{equation*}

We denote this multiscale Bloch mapping by $m_{S}=(m_{l,S})_{l \in S}$, which approximates the true Bloch mapping $m=(m_{l})_{l=1}^{L}$.

\subsection{Multiscale Bloch mapping derivative}

From (3), it follows that
\begin{equation*}
    m^{\prime}_{l+N}=\left(\partial_{T_1} m_{l+N}, \partial_{T_2} m_{l+N}, \partial_{\omega} m_{l+N}  \right)
\end{equation*}
using the notation that $\partial_{\chi} m=\partial m/\partial \chi$ for $\chi$ denoting one of the variables $(T_1,T_2,\omega)$. In particular,
\begin{multline*}
    \partial_\chi m_{l+N}=\partial_\chi A^{N} m_{l} + A^{N} \partial_\chi m_{l} + \sum_{k=0}^{N-1} \partial_\chi A^{k} b + A^{k} \partial_\chi b.
\end{multline*}
Using the representation of $A$ in (4), we have
\begin{equation*}
    \partial_\chi A^{k}= \sum_{i=1}^{3} k \lambda_{i}^{k-1} \partial_\chi \lambda_{i} U_{i}+\lambda_{i}^k \partial_\chi U_{i},
\end{equation*}
where
\begin{equation*}
\partial_\chi U_{i}=\left(\partial_\chi Z G_{i} - U_{i} \partial_\chi Z\right)Z^{-1}.
\end{equation*}

\subsection{Objective functions}
Analogous to (2), we define $Q_{l,S} (\rho,u) = P_{l}  \mathcal{F} (\rho T_{x,y} m_{l,S}(u))$ for the multiscale version of the forward operator. To reconstruct the tissue parameters $x=(\rho,u)$, we minimise the multiscale objective function
\begin{equation}
F_{S}(x)=\frac{1}{2 |S|} \sum_{l \in S} \| Q_{l,S}(x)-y_{l}\|^{2},
\end{equation}
where $\| \cdot \|$ denotes the usual complex norm. This approximates the true objective function $F$, given by $F=F_{S_1(1)}$.

The corresponding gradient of (5) is calculated as
\begin{equation*}
    \nabla F_{S}(x)=\frac{1}{|S|} \sum_{l \in S}  \left(Q_{l,S}^{\prime}(x)\right)^{*} \left[ Q_{l,S}(x)-y_{l} \right],
\end{equation*}
where ${\,}^{*}$ denotes the conjugate transpose and the Jacobian of the multiscale forward operator is a linear operator that acts on a vector $h=(h_{\rho},h_{u})^{T}$ where $h_{u}=(h_{T_1}, h_{T_2}, h_{\omega})^{T}$, as follows
\begin{equation*}
Q^{\prime}_{l,S}(x) h = P_{l} \mathcal{F} \, T_{x,y}(h_{\rho} m_{l,S}(u)+\rho m_{l,S}^{\prime}(u) h_{u}).
\end{equation*}

\section{Algorithms}
\label{sec:methods}

In this chapter, we discuss algorithms to minimise our objective functions. 

\subsection{Projected coordinate descent with backtracking}

We will need the following backtracking condition
\begin{equation}
    F_{S}(x^{k+1}) \leq F_{S}(x^{k}) + \langle \partial_i F_{S}(x^{k}), \Delta_{i}^{k+1} \rangle +  \frac{\|\Delta_{i}^{k+1}\|^2}{2 \tau_{i}^{k}}
\end{equation}
and a projection $\pi_{i}(x)=\argmin_{y \in [\theta, \infty)} \|y-x\|$ (e.g. $\theta=1$) when applied to tissue parameter maps that take non-negative values (e.g. $\rho$, $T_1$, $T_2$) and $\pi_{i}(x)=x$ otherwise (e.g. $\omega$). Note that $\partial_{i} F_{S}(x^{k})$ corresponds to the $i$-th component of $\nabla F_{S}(x^{k})$, $\langle \cdot,\cdot \rangle$ denotes the standard Euclidean inner product and $\Delta_{i}^{k+1}=x_{i}^{k+1}-x_{i}^{k}$.

\begin{algorithm}
\SetKwInput{Inputs}{input}
\Inputs{$x^{0} \in \mathds{R}^{n \times m \times I}$, $\tau^{0} \in \mathds{R}^I$, $N, K \in \mathds{R}$}
\For{$k=0$ \KwTo $K-1$}{
$x^{k+1}=x^{k}$\;
$S=S_N(\operatorname{randi}(N))$\;
\For{$i=1$ \KwTo $I$}{
\For{$c=1$ \KwTo $C$}{
$x_{i}^{k+1} = \pi_{i} \left(x_{i}^{k}-\tau_{i}^{k} \partial_{i} F_{S}(x^{k}) \right)$\;
\eIf{(6) is satisfied}{
$x_{i}^{k} =x_{i}^{k+1}$\;
$\tau_{i}^{k+1} =\overline{\eta}  \tau_{i}^{k}$\;
Break inner for loop\;
}{
$\tau_{i}^{k}=\underline{\eta} \tau_{i}^{k}$\;
\If{$c=C$}{
$\tau_{i}^{k+1}=\tau_{i}^{k}$\;
}
}
}
}
}
\SetKwInput{Output}{output}
\Output{$\left(x^{K}, \tau^{K}\right)$}
\caption{Projected coordinate descent with backtracking (PCDB)}
\end{algorithm}

The inputs $x^{0}$, $\tau^{0}$, $N$ and $K$ correspond to the initial guess, the initial step size, the grid increment and the number of iterations, respectively. Backtracking involves updating the step size component-wise at each iteration, at most $C$ times (e.g. $C=50$), via $\overline{\eta}$ and $\underline{\eta}$ (e.g. $\overline{\eta}=1.2$, $\underline{\eta}=0.75$) to ensure monotonic descent of the objective function. Note that $\operatorname{randi}(N)$ represents choosing an integer uniformly at random from the set $\{1,...,N\}$.

\subsection{Coarse-to-fine}
We let $N$ be a vector of natural numbers of length $J$, with entries that are strictly decreasing, e.g. $N=(25,10,5,1)$. The entries in $N$ specify the increments of uniform grids that are then used in algorithm 1. Moreover, we let $K$ be another vector of natural numbers of length $J$, whose entries correspond to the number of iterations we use on each grid, e.g. $K=(250,200,100,950)$. 

\begin{algorithm}
\SetKwInOut{Inputs}{input}
\Inputs{$x \in \mathds{R}^{n \times m \times I}$, $\tau \in \mathds{R}^I$, $N, K \in \mathds{R}^J$}
\For{$j=1$ \KwTo $J$}{
$\left(x,\tau\right)$ = PCDB($x$, $\tau$, $N(j)$, $K(j)$)\;
}
\SetKwInOut{Output}{output}
\Output{$x$}
\caption{Coarse-to-fine (C2F)}
\end{algorithm}

\section{Numerical results}
\label{sec:results}

\begin{table}[b]
\caption{PSNR (dB) of reconstructed tissue parameter maps.}
\centering
 \begin{tabular}{ c  c  c  c  c }
  Method & $\rho$ & $T_1$ & $T_2$ & $\omega$ \\
 \hline
 BLIP & $28.69$ & $31.08$ & $32.00$ & $36.01$ \\

 FINE & $32.25$ & $21.33$ & $28.59$ & $20.44$  \\

 C2F & $37.63$ & $22.86$ & $28.32$ & $21.87$  \\

 BLIP+FINE & $50.80$ & $35.56$ & $33.02$ & $64.35$  \\

 BLIP+C2F & $50.97$ & $35.69$ & $33.24$ & $64.01$ 

\end{tabular}
\end{table}

We now present two experiments comparing the C2F algorithm with an analogous approach where no multiscaling is used, i.e. algorithm 1 with $N=1$ (FINE). One experiment uses a constant initial guess, whereas the other uses an initial guess generated by BLIP \cite{davies2014compressed, golbabaee2017cover} with a coarse dictionary simulated from the discretisation $T_1=\{0.5,1,...,6\}$ (s), $T_2=\{0.05,0.1,...,0.6\}$ (s) and $\omega=\{-50,-40,...,50\}$ (Hz). Our tests retrospectively simulate MRF data from the ground truth quantitative parameter maps of a healthy volunteer's brain, provided in \cite{ma2013magnetic}. Images have a spatial resolution of $128$$\times$$128$ pixels. We consider a randomised Cartesian subsampling operator [3], in particular, we use an $n$-multishot Echo Planar Imaging (EPI) scheme with random offsets on each frame. A subsampling rate of $1/8$ is used throughout the tests presented.

Fig. 2 displays the true objective function $F$ evaluated at each iterate $x^k$ in the two experiments. Note that the objective graphs corresponding to the C2F algorithm are scaled (in the iteration axis) by $1/N$, where $N$ is the increment of the uniform grid, with vertical partitions indicating positions of refinement. This gives a measure of the relative computational cost between the two methods and is chosen so that the total number of gradients calculated in both algorithms is the same. Fig. 3 then displays the reconstructed tissue parameter maps (odd columns) with the corresponding error maps (even columns).

We observe that C2F successfully minimises the true objective function at a faster rate than FINE in both experiments. This is also reflected in tables 1 and 2, which give the Peak Signal-to-Noise Ratio (PSNR) and the Mean Absolute Percentage Error (MAPE) of final reconstructions parameter-wise. We see that C2F results in consistently lower percentage errors and generally higher PSNR values, with the latter indicating higher image quality. Note that due to the cyclic nature of $\omega$, pixels that converge to the true value $\pm$ the period ($1/T_R=100$ Hz) are considered accurate since they produce the same Bloch responses (see (1)). MAPE is not given for $\omega$ because of the $0$-values in the ground truth. Whilst the accuracy of BLIP is tied to its dictionary size (and is inaccurate here because of this), one can still use BLIP to generate the initial guess for C2F, as we have done, to get superior accuracy without sacrificing memory. 

\begin{figure}[t]

\begin{minipage}[b]{1.0\linewidth}
  \centering
  \centerline{\includegraphics[width=8.5cm,trim={0cm 0cm 0cm 0cm},clip]{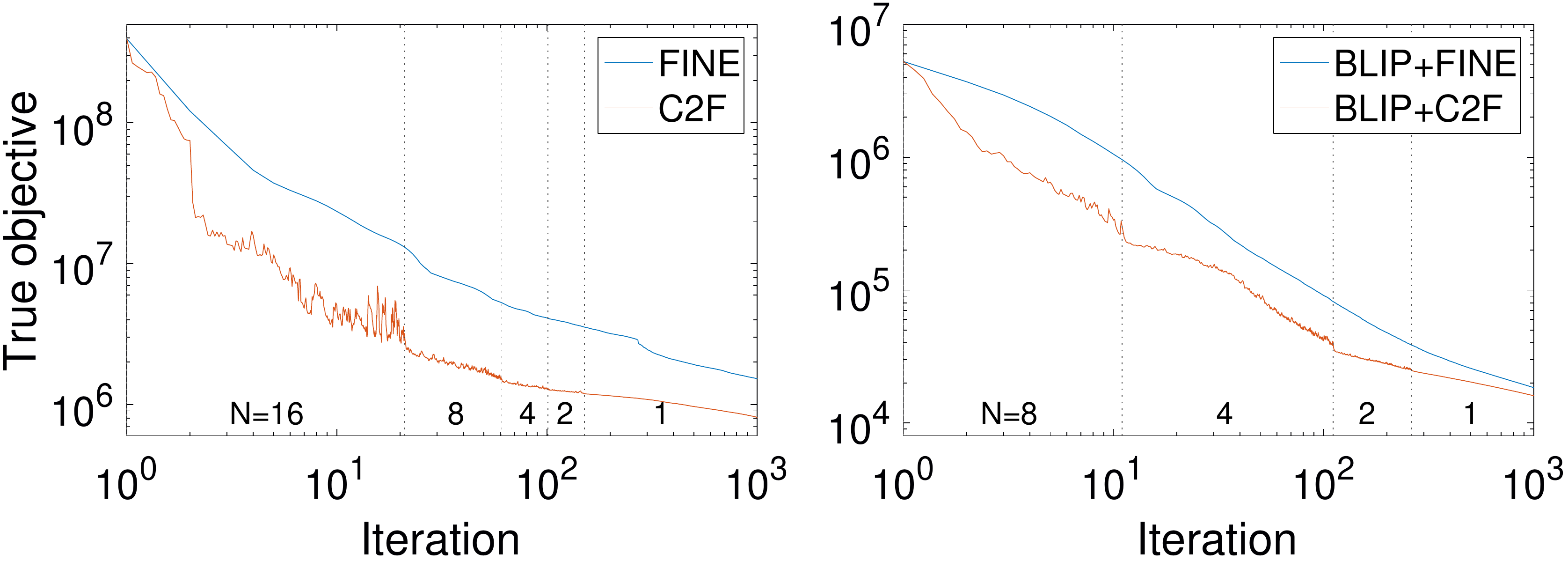}}
%  \vspace{2.0cm}
\end{minipage}

\caption{True objective comparison of FINE with C2F (left), and BLIP+FINE with BLIP+C2F (right). On the left, $x^0$ consists of four constant images taking values $(\rho,T_1,T_2,\omega)=(0.42,2,0.2,0)$ (see fig. 3 for units), whereas on the right $x^0$ is generated by BLIP with a coarse dictionary (see fig. 3, row 1). In both, the initial step size is $\tau^0=(0.1,1,0.1,10^{-8})$. The remaining inputs are $N=1$, $K=1000$ for FINE and BLIP+FINE, $N=(16,8,4,2,1)$, $K=(320,320,160,100,850)$ for C2F, and $N=(8,4,2,1)$, $K=(80,400,300,740)$ for BLIP+C2F.}
\label{fig:SS}
\end{figure}

\begin{table}[t]
\caption{MAPE ($\%$) of reconstructed tissue parameter maps.}
\centering
 \begin{tabular}{ c  c  c  c  }
 Method & $\rho$ & $T_1$ & $T_2$ \\
 \hline
 BLIP & $14.97$ & $\phantom{0}45.62$ & $23.73$ \\

 FINE & $11.87$ & $162.79$ & $30.01$  \\

 C2F & $\phantom{0}9.67$ & $122.33$ & $15.29$  \\

 BLIP+FINE & $\phantom{0}2.72$ & $\phantom{00}2.69$ & $\phantom{0}3.85$  \\

 BLIP+C2F & $\phantom{0}2.56$ & $\phantom{00}2.51$ & $\phantom{0}3.54$ 

\end{tabular}
\end{table}

\begin{figure*}[t]
\begin{center}
\begin{tikzpicture}
\draw (0cm, 0cm) node[anchor=south] {\includegraphics[width=2.1cm, height=2.1cm, clip, trim={1cm 1cm 1cm 1cm}]{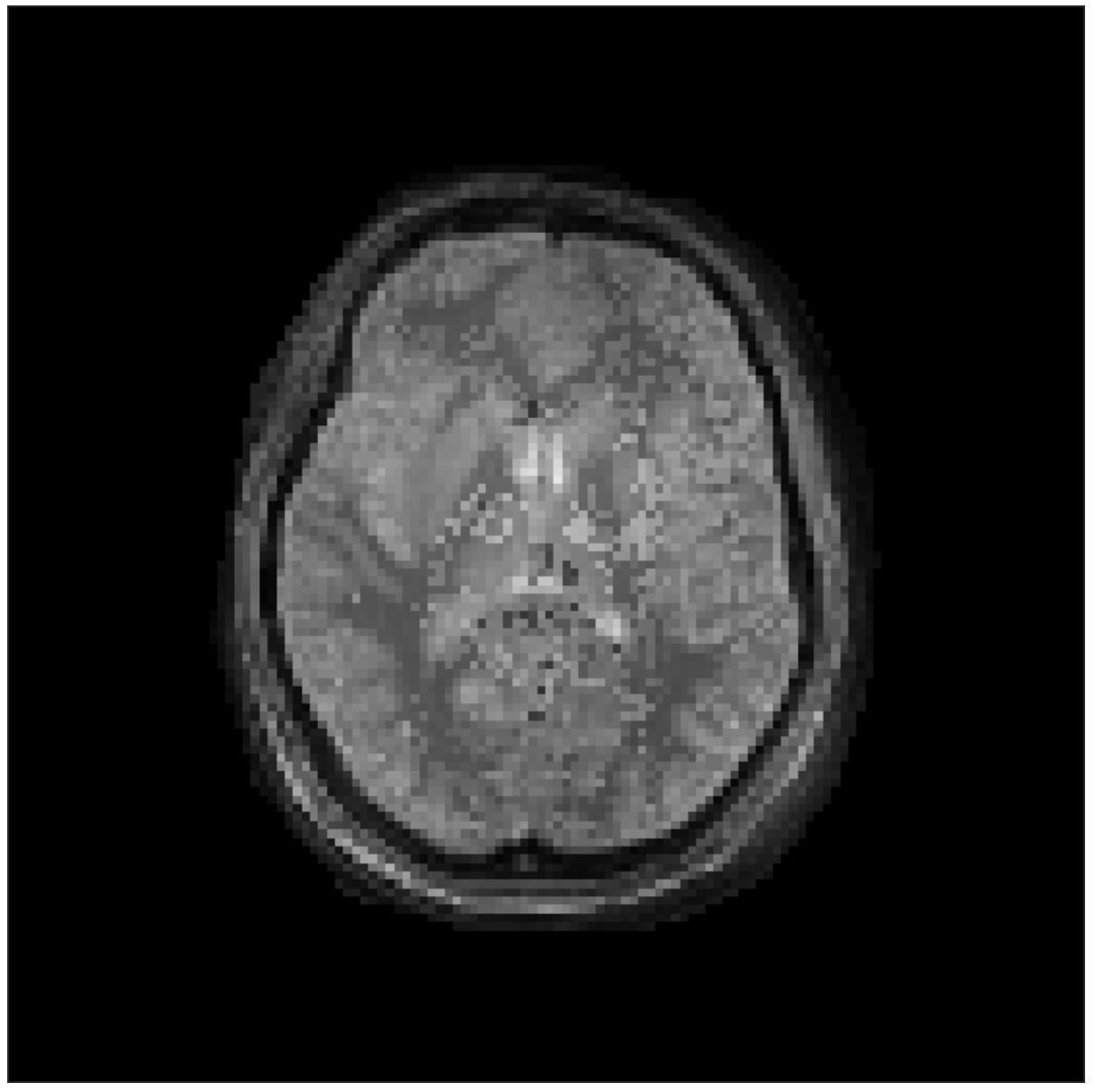}};

\draw (2.1cm, 0cm) node[anchor=south] {\includegraphics[width=2.1cm, height=2.1cm, clip, trim={1cm 1cm 1cm 1cm}]{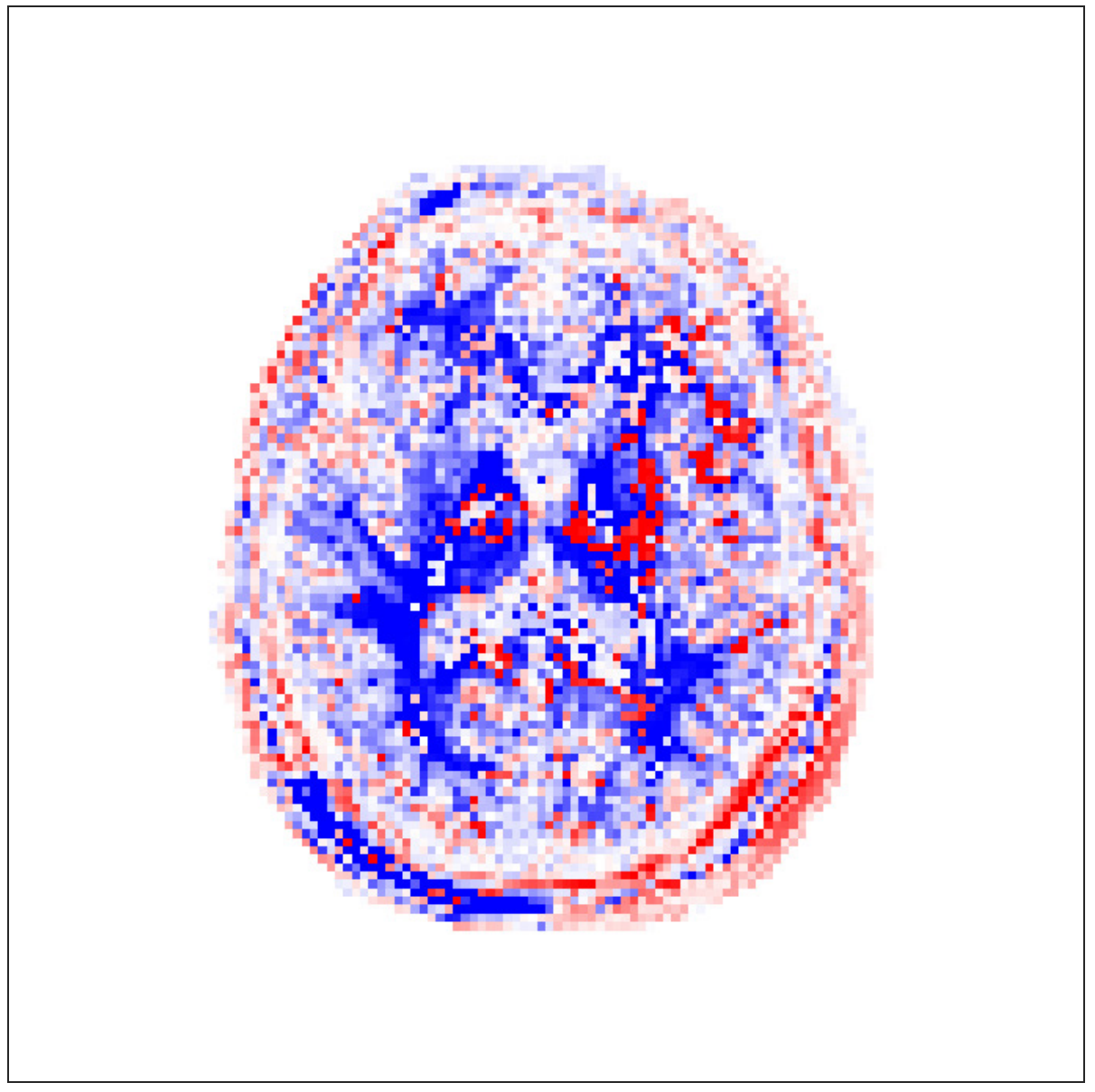}};

\draw (4.2cm, 0cm) node[anchor=south] {\includegraphics[width=2.1cm, height=2.1cm, clip, trim={1cm 1cm 1cm 1cm}]{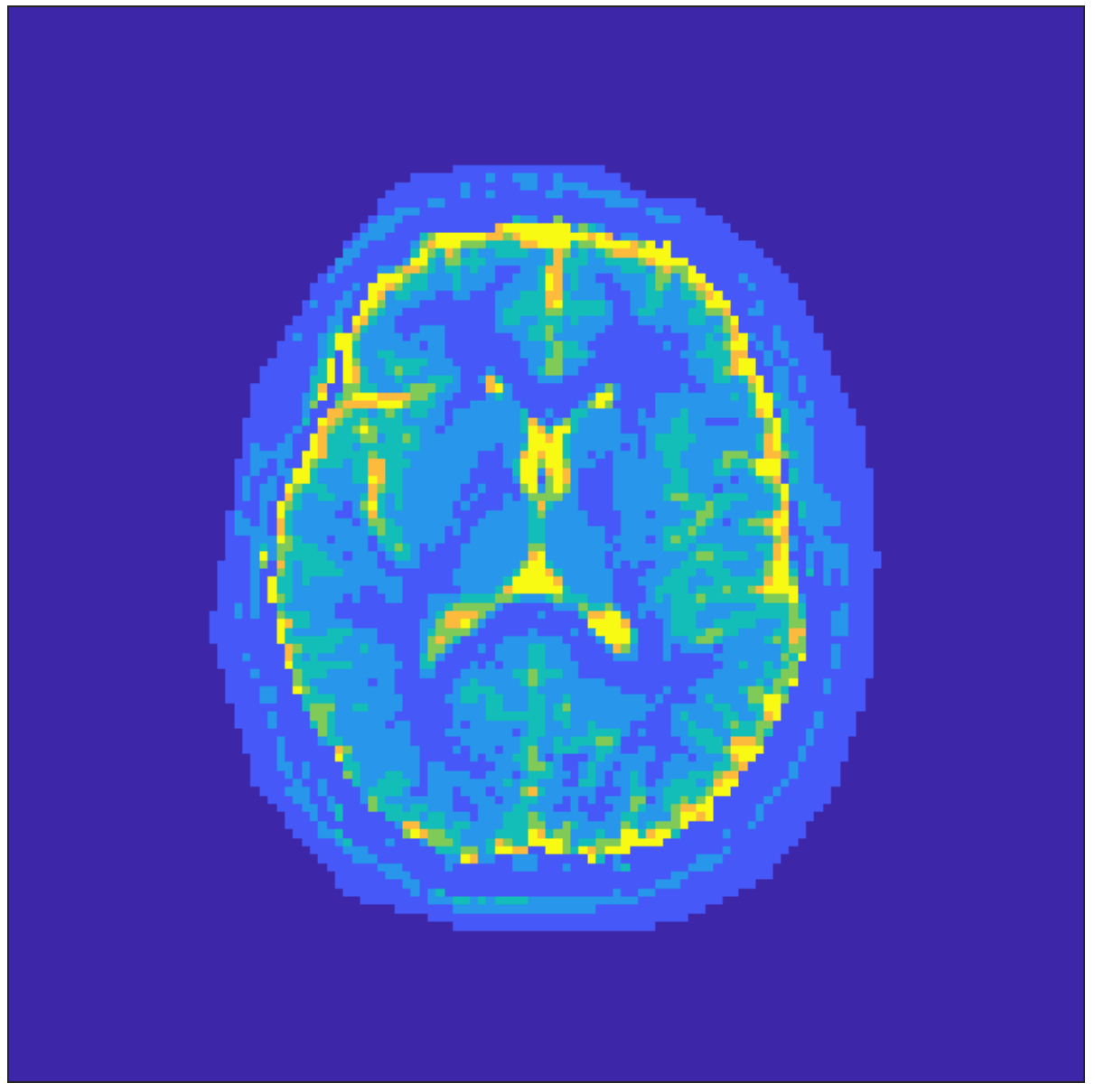}};

\draw (6.3cm, 0cm) node[anchor=south] {\includegraphics[width=2.1cm, height=2.1cm, clip, trim={1cm 1cm 1cm 1cm}]{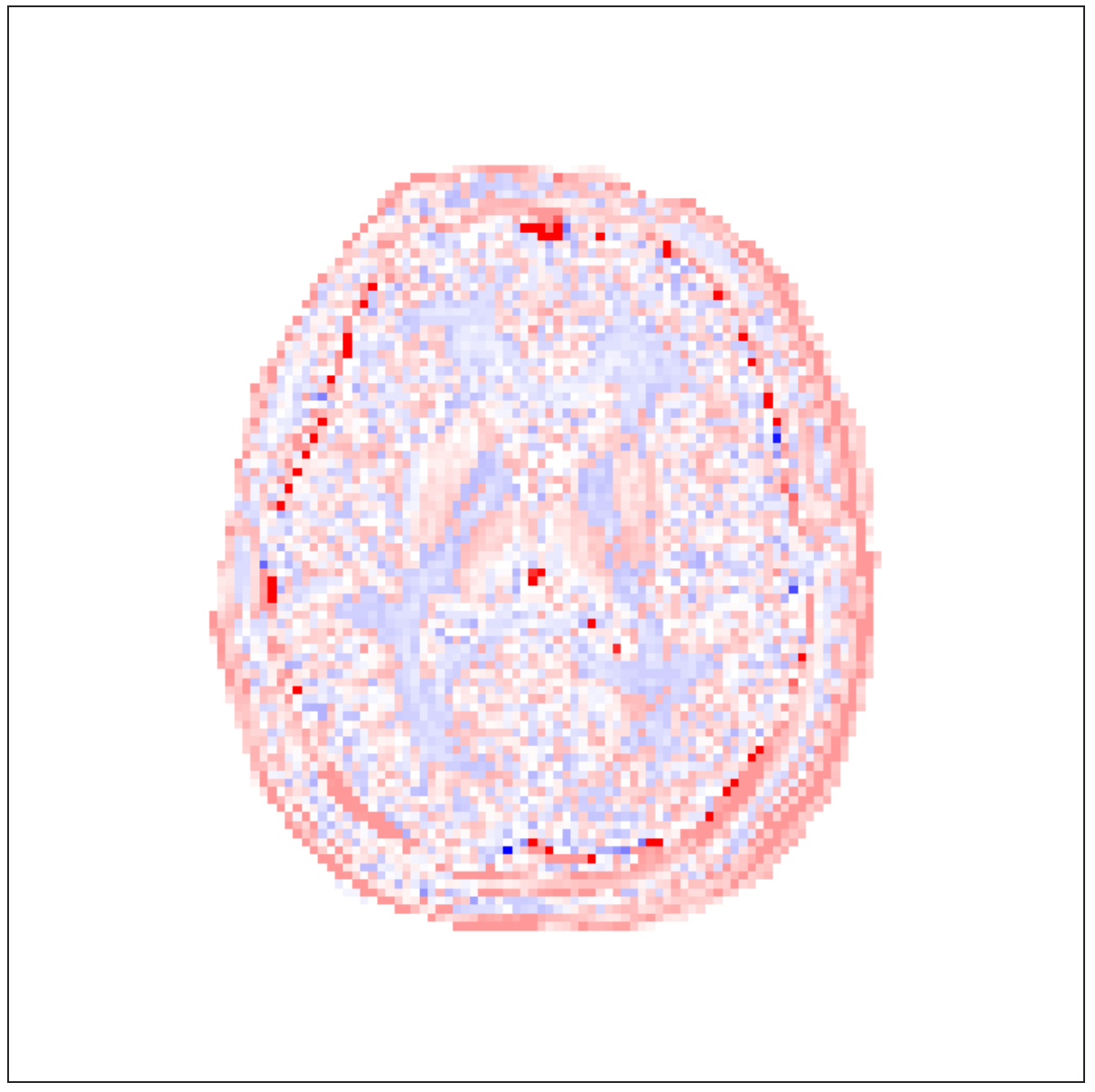}};

\draw (8.4cm, 0cm) node[anchor=south] {\includegraphics[width=2.1cm, height=2.1cm, clip, trim={1cm 1cm 1cm 1cm}]{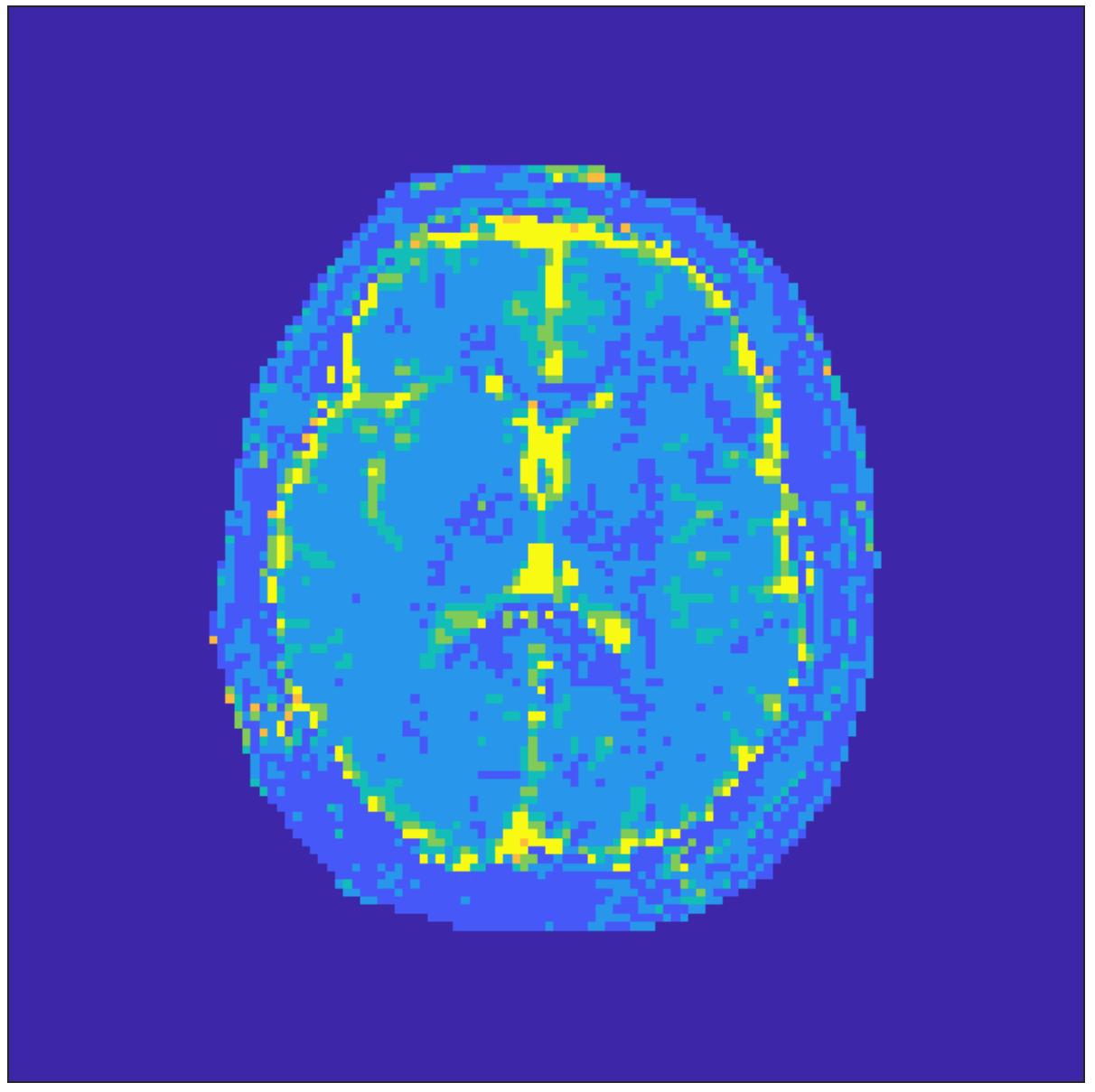}};

\draw (10.5cm, 0cm) node[anchor=south] {\includegraphics[width=2.1cm, height=2.1cm, clip, trim={1cm 1cm 1cm 1cm}]{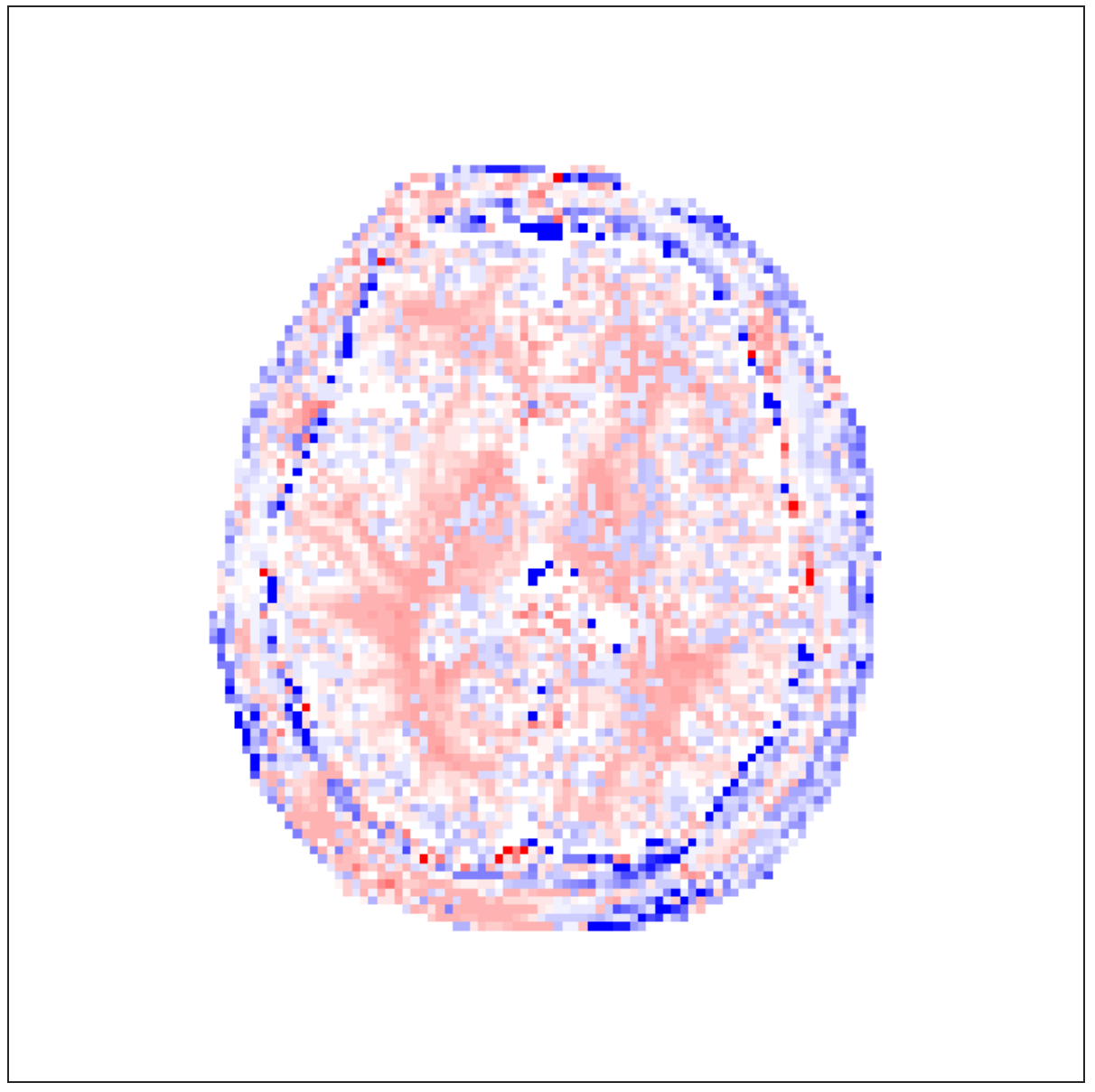}};

\draw (12.6cm, 0cm) node[anchor=south] {\includegraphics[width=2.1cm, height=2.1cm, clip, trim={1cm 1cm 1cm 1cm}]{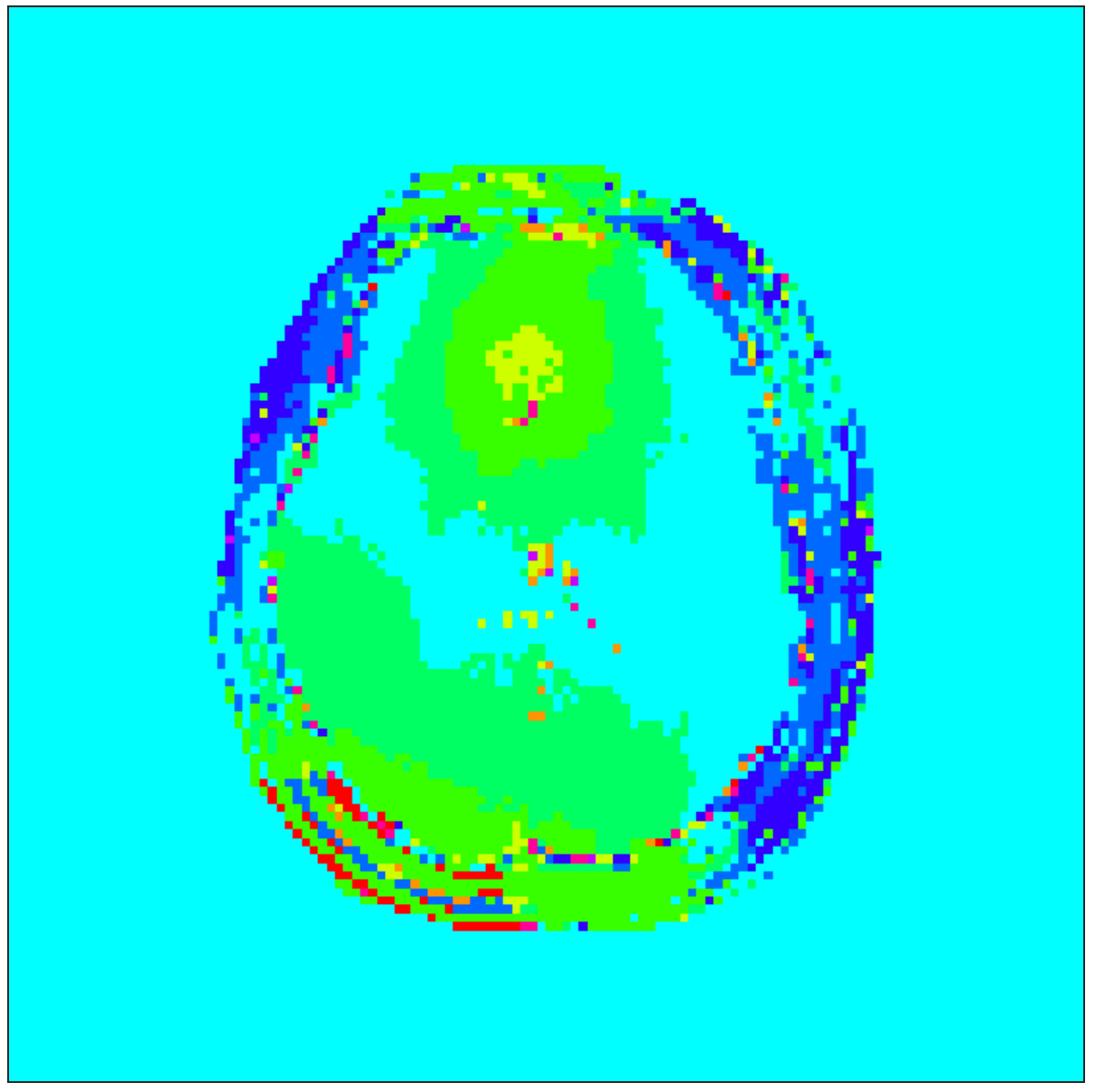}};

\draw (14.7cm, 0cm) node[anchor=south] {\includegraphics[width=2.1cm, height=2.1cm, clip, trim={1cm 1cm 1cm 1cm}]{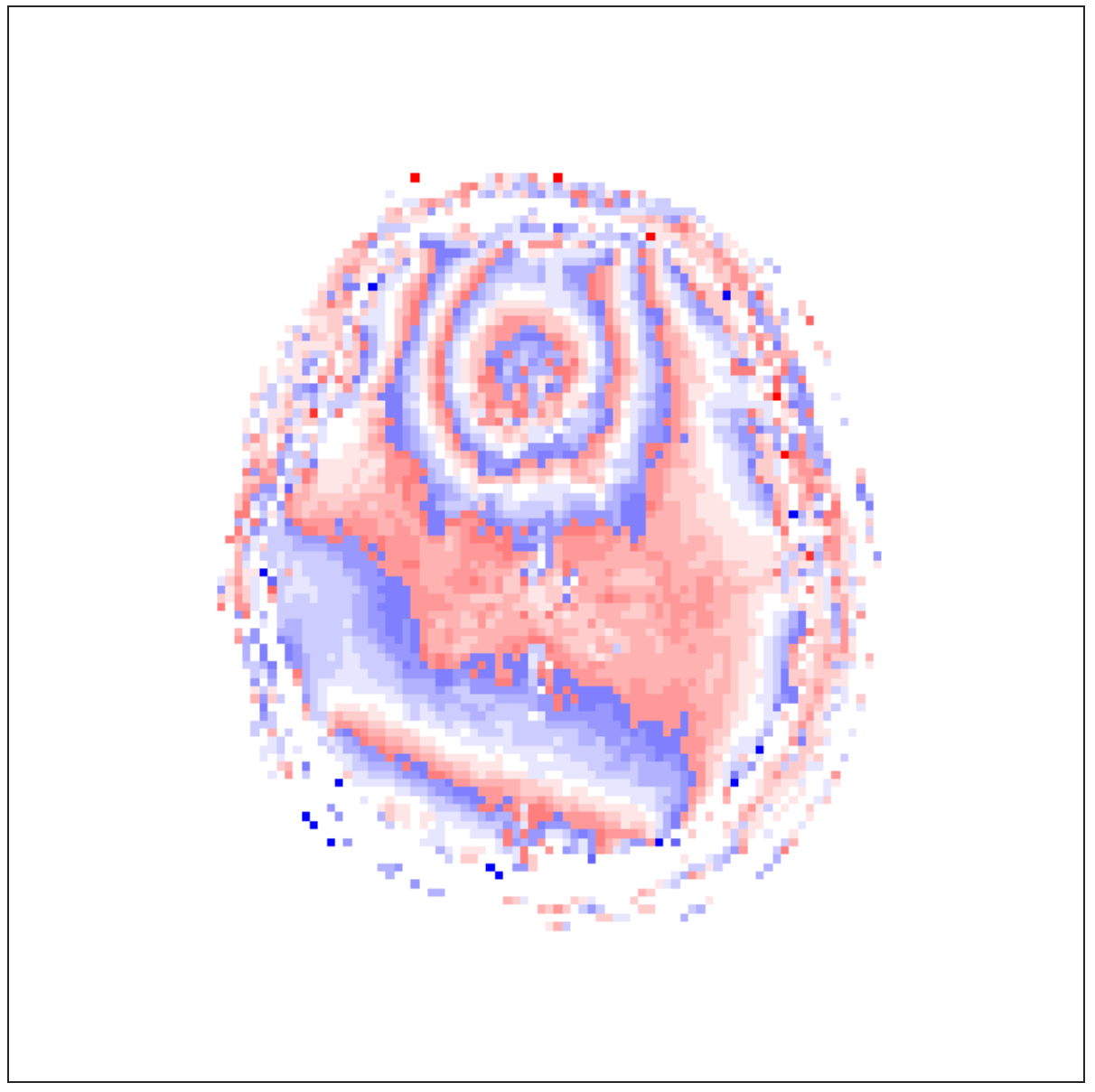}};
%%% new row
\draw (0cm, -2.15cm) node[anchor=south] {\includegraphics[width=2.1cm, height=2.1cm, clip, trim={1cm 1cm 1cm 1cm}]{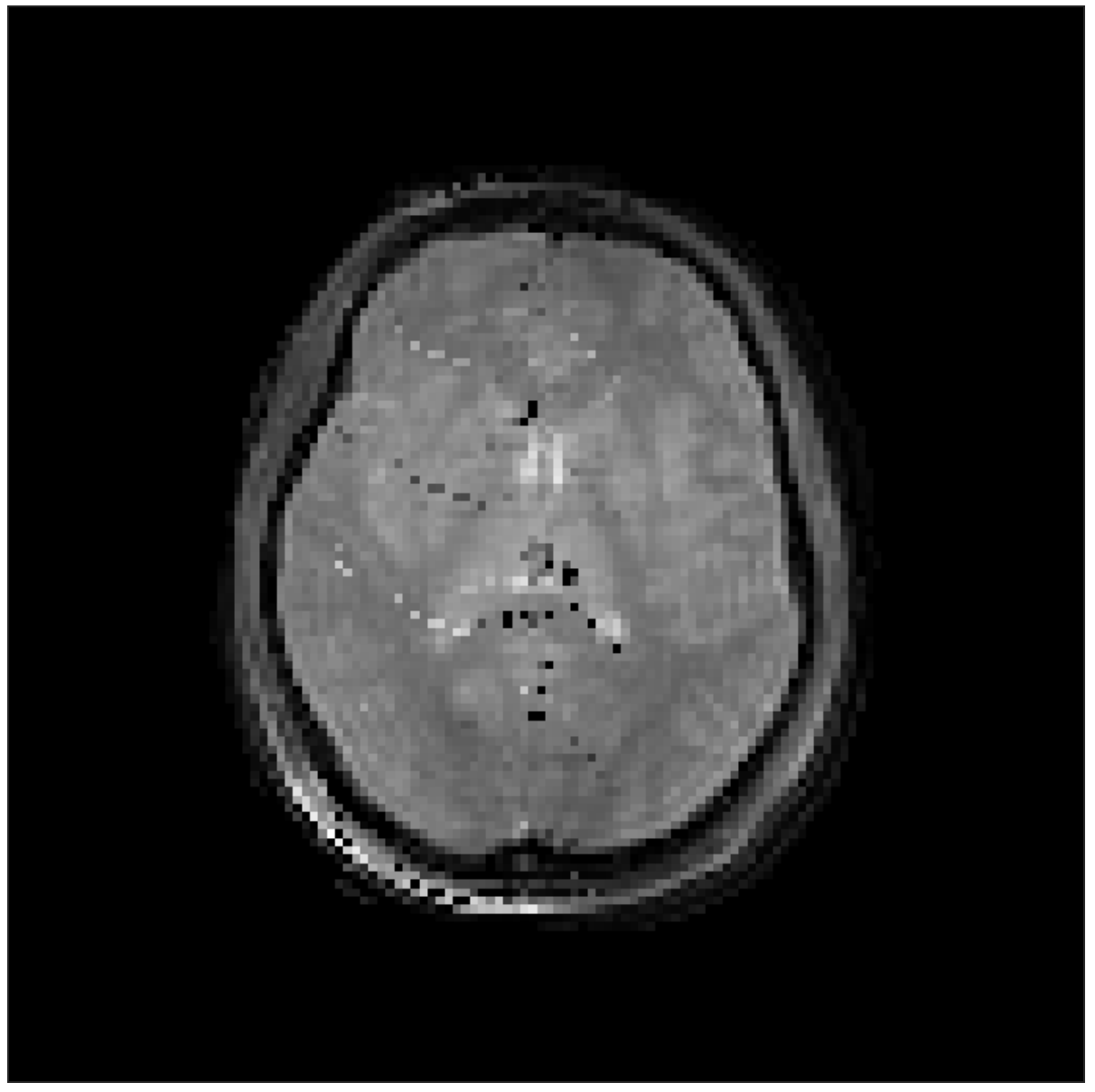}};

\draw (2.1cm, -2.15cm) node[anchor=south] {\includegraphics[width=2.1cm, height=2.1cm, clip, trim={1cm 1cm 1cm 1cm}]{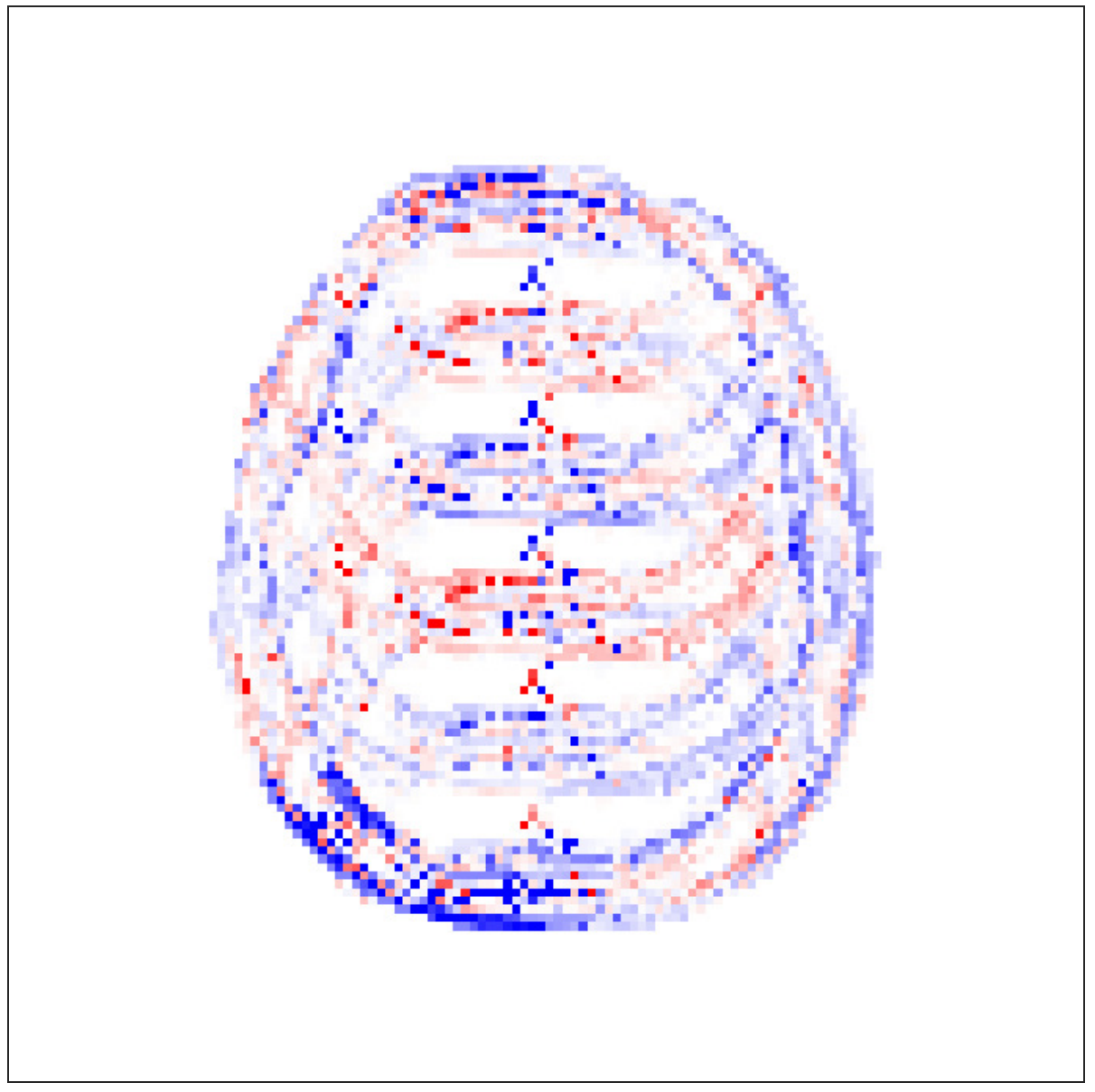}};

\draw (4.2cm, -2.15cm) node[anchor=south] {\includegraphics[width=2.1cm, height=2.1cm, clip, trim={1cm 1cm 1cm 1cm}]{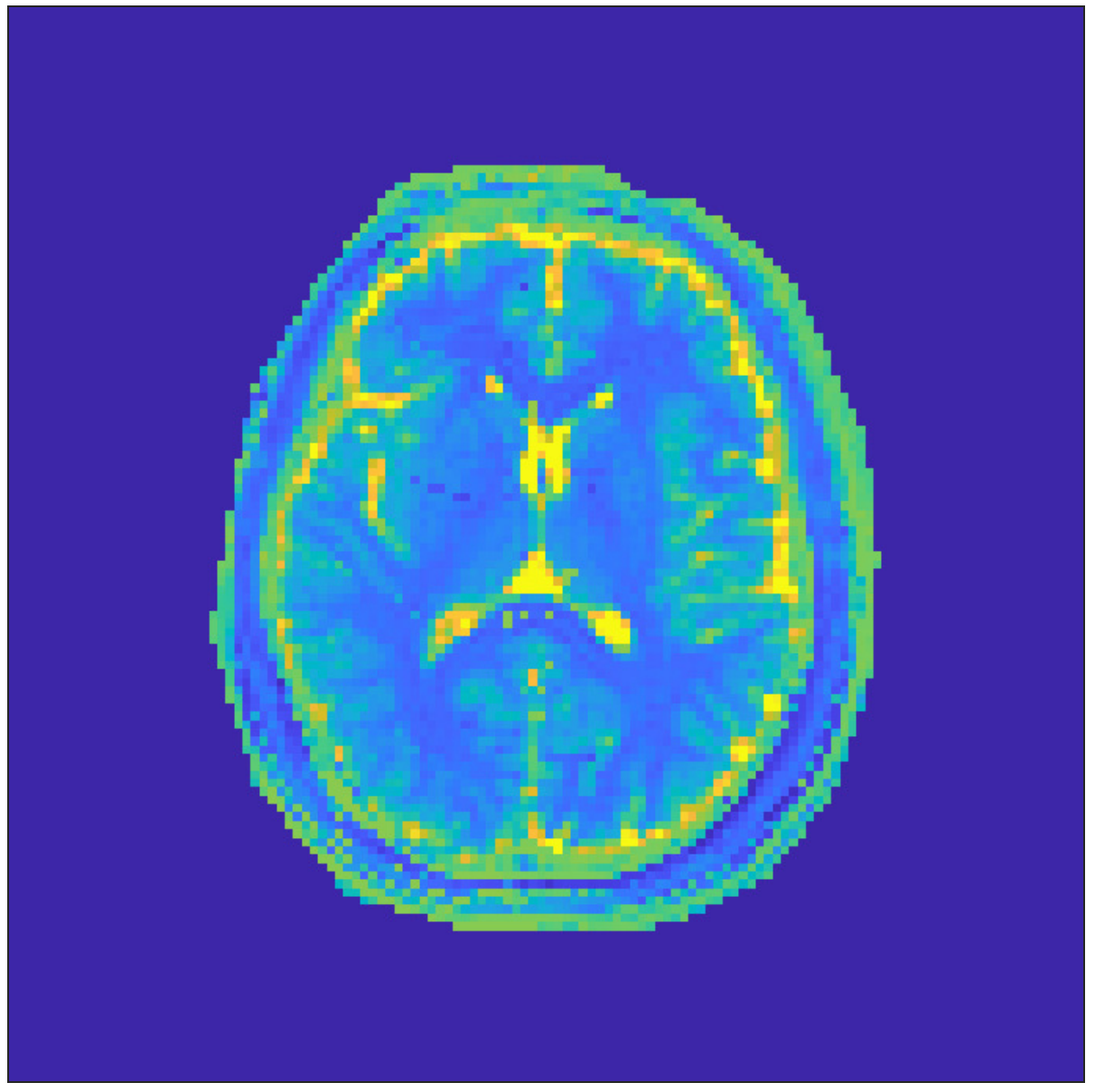}};

\draw (6.3cm, -2.15cm) node[anchor=south] {\includegraphics[width=2.1cm, height=2.1cm, clip, trim={1cm 1cm 1cm 1cm}]{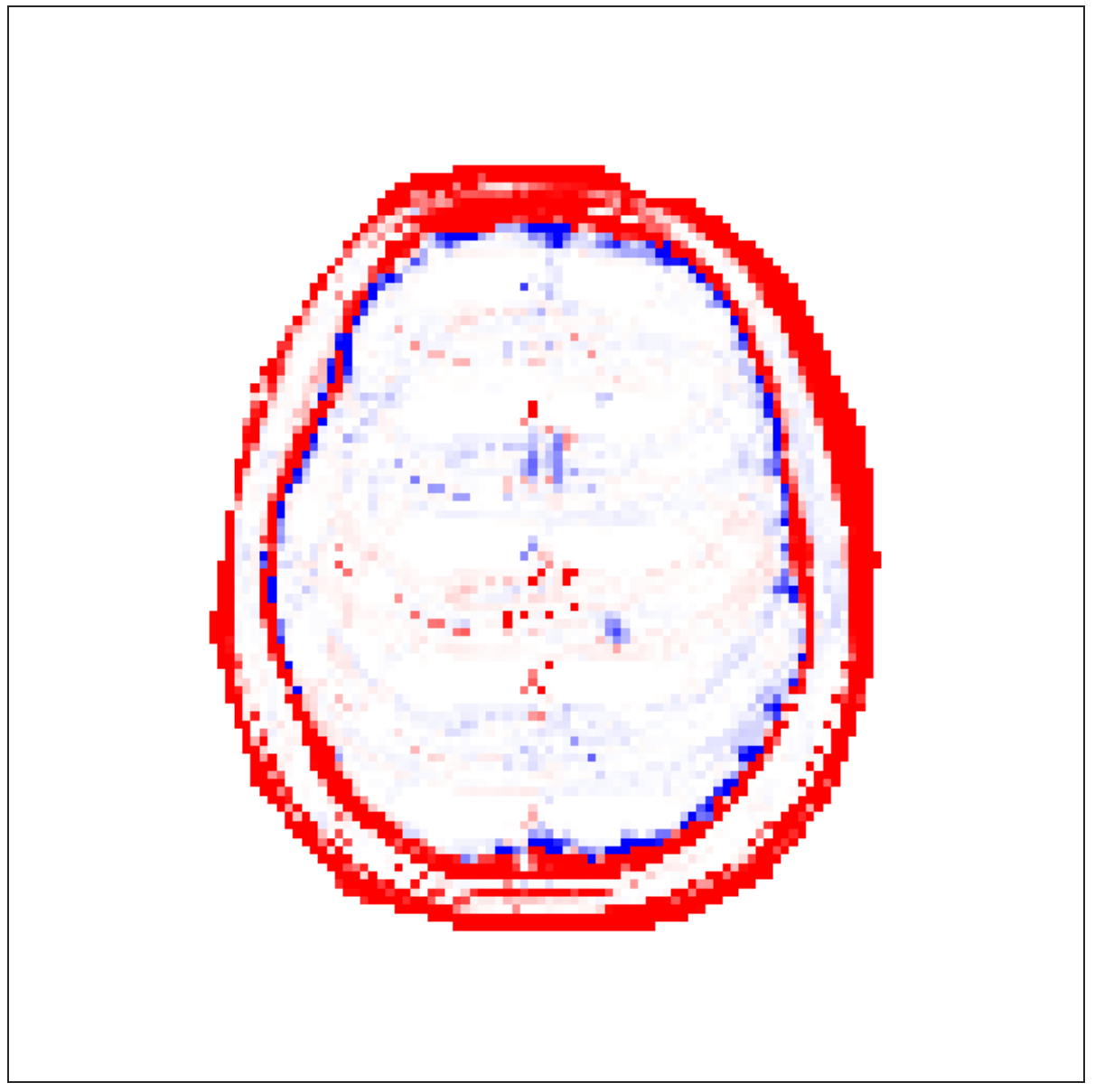}};

\draw (8.4cm, -2.15cm) node[anchor=south] {\includegraphics[width=2.1cm, height=2.1cm, clip, trim={1cm 1cm 1cm 1cm}]{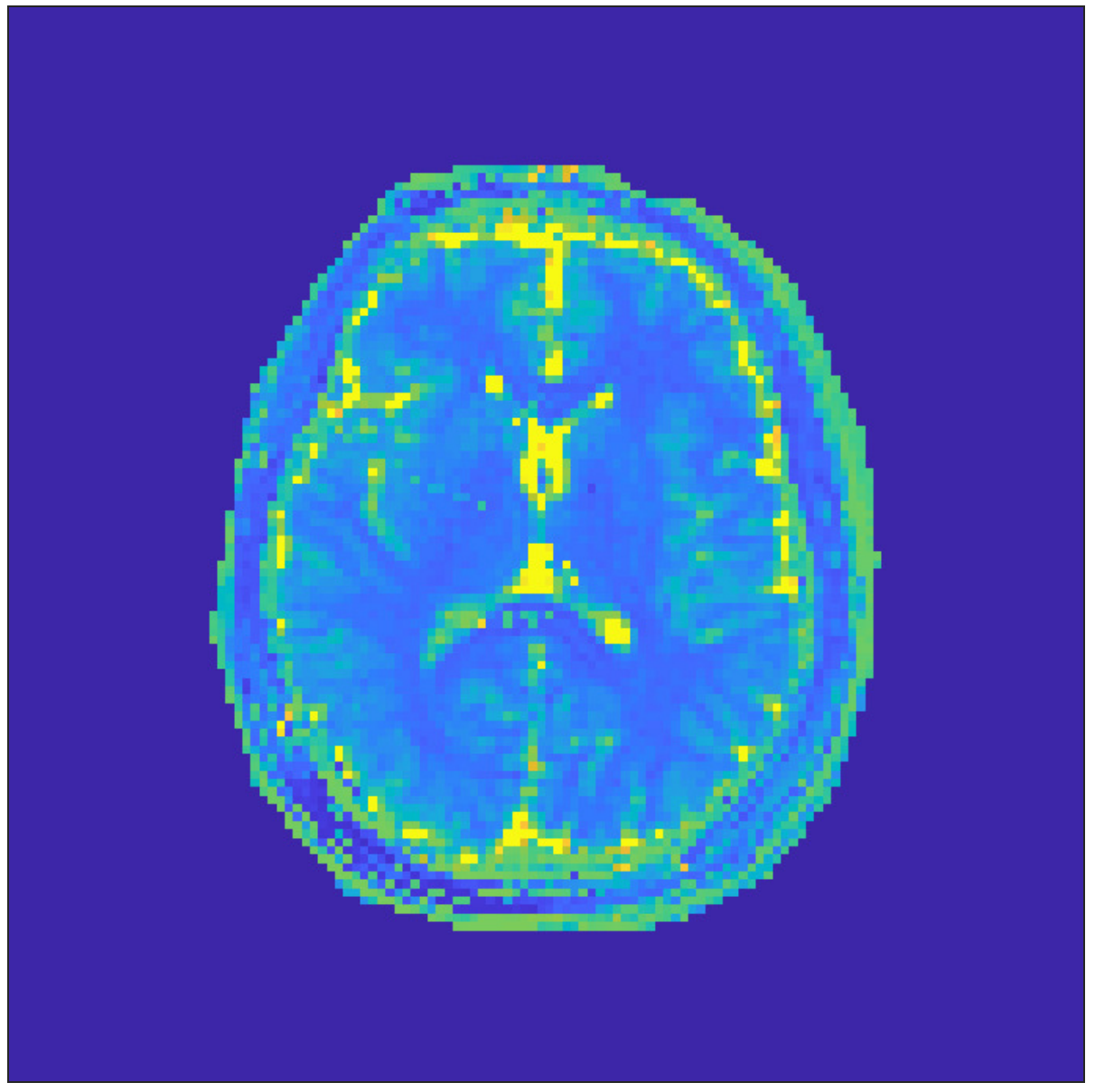}};

\draw (10.5cm, -2.15cm) node[anchor=south] {\includegraphics[width=2.1cm, height=2.1cm, clip, trim={1cm 1cm 1cm 1cm}]{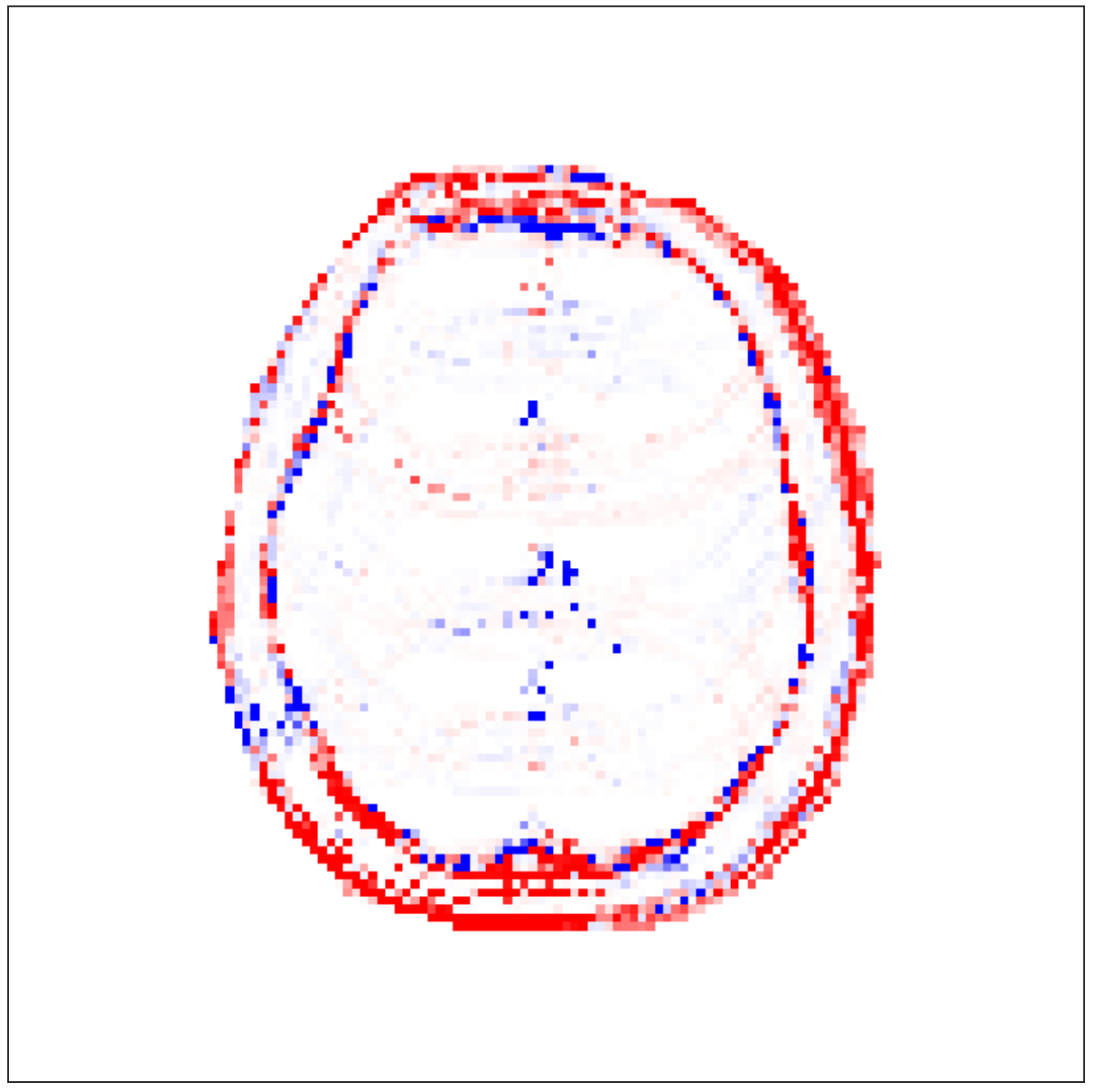}};

\draw (12.6cm, -2.15cm) node[anchor=south] {\includegraphics[width=2.1cm, height=2.1cm, clip, trim={1cm 1cm 1cm 1cm}]{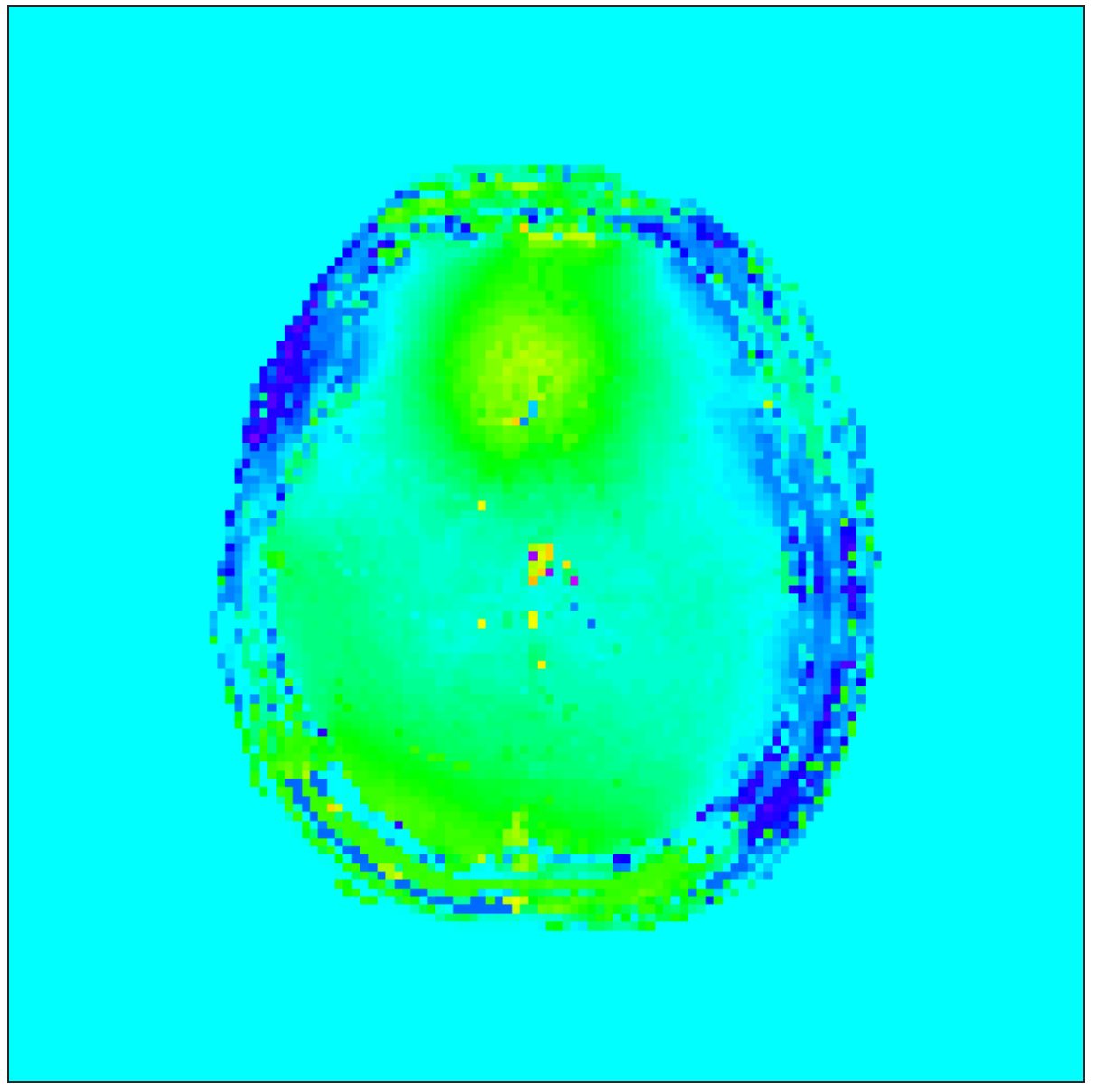}};

\draw (14.7cm, -2.15cm) node[anchor=south] {\includegraphics[width=2.1cm, height=2.1cm, clip, trim={1cm 1cm 1cm 1cm}]{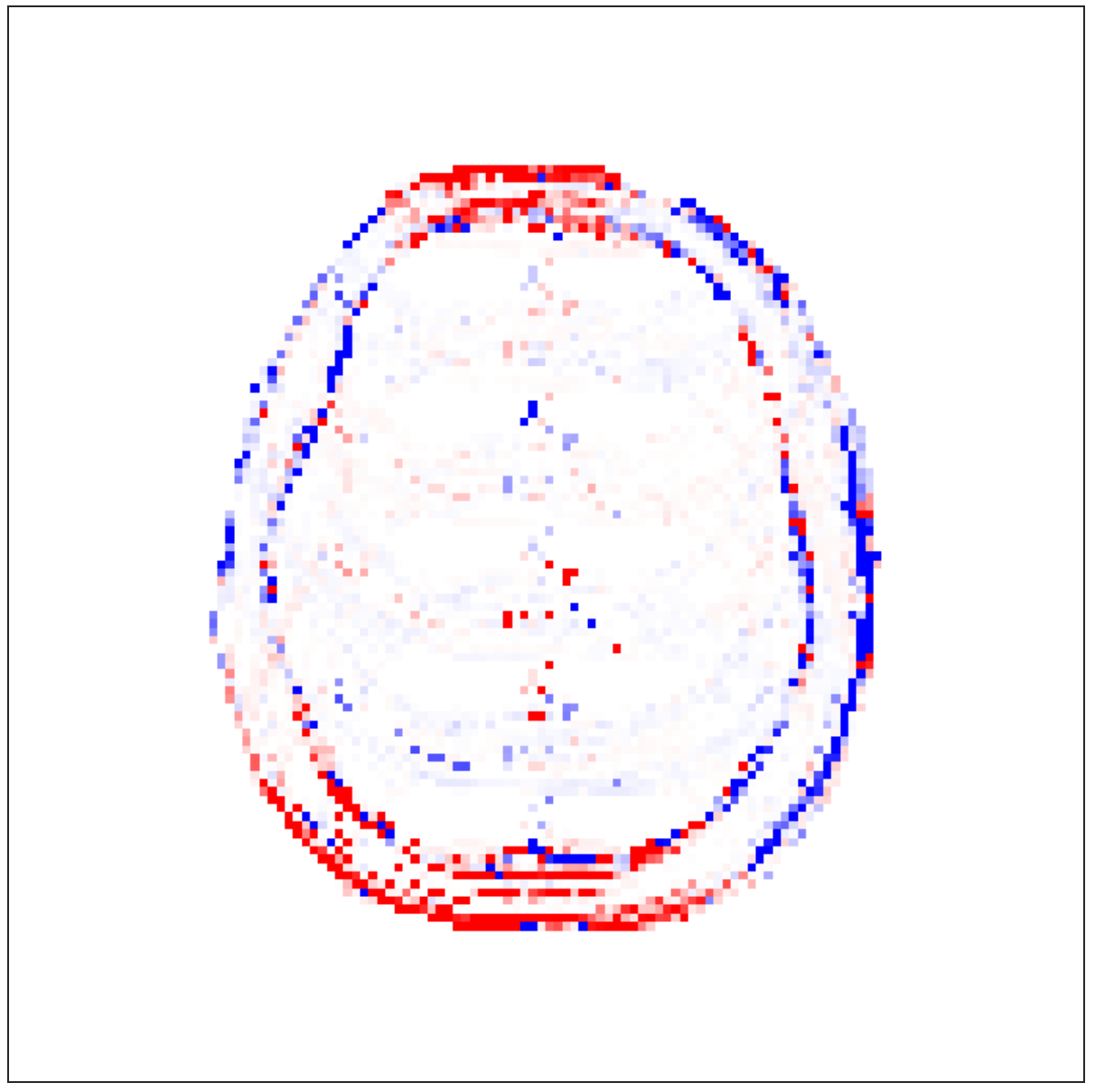}};
%%% new row
\draw (0cm, -4.3cm) node[anchor=south] {\includegraphics[width=2.1cm, height=2.1cm, clip, trim={1cm 1cm 1cm 1cm}]{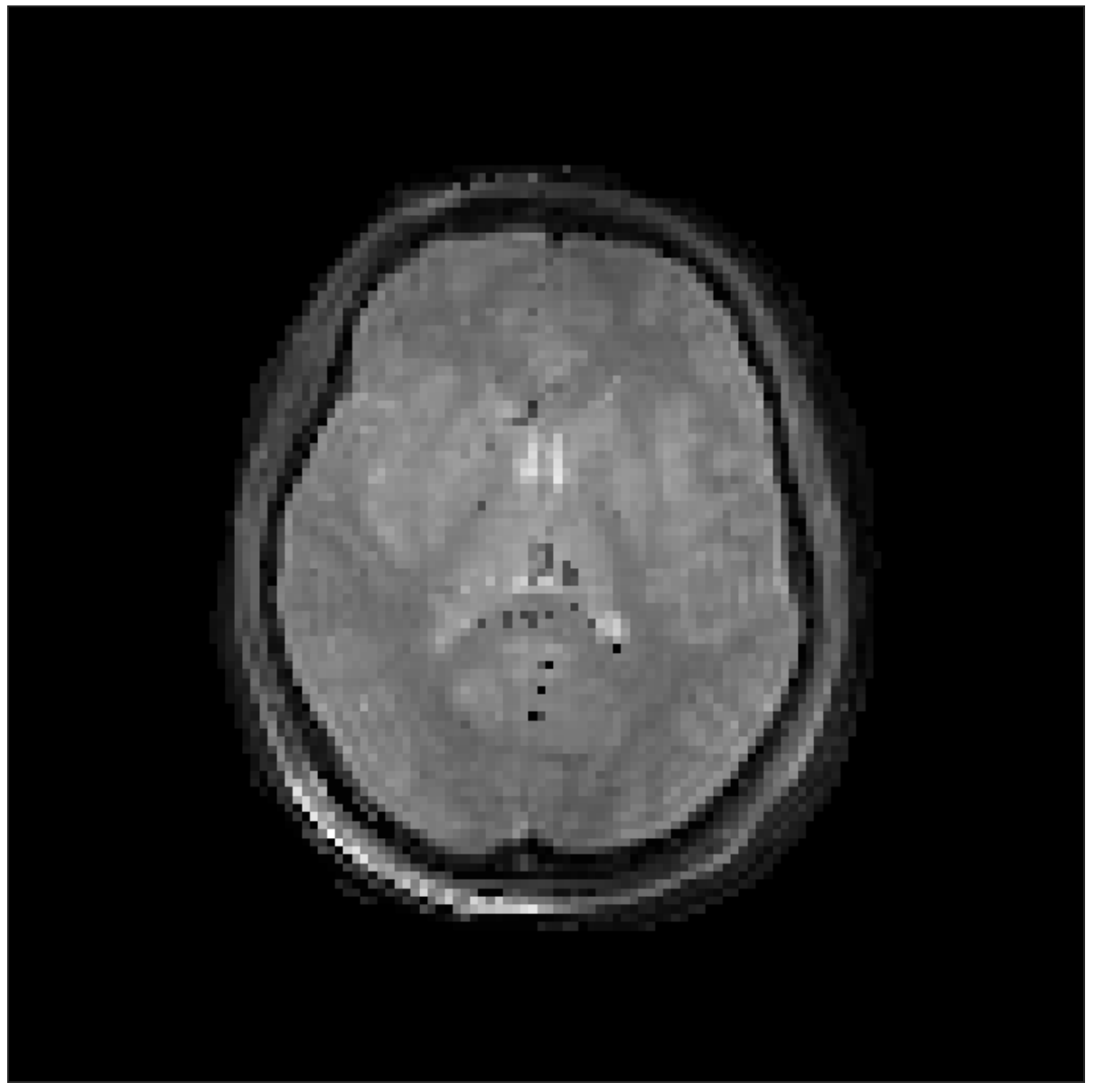}};

\draw (2.1cm, -4.3cm) node[anchor=south] {\includegraphics[width=2.1cm, height=2.1cm, clip, trim={1cm 1cm 1cm 1cm}]{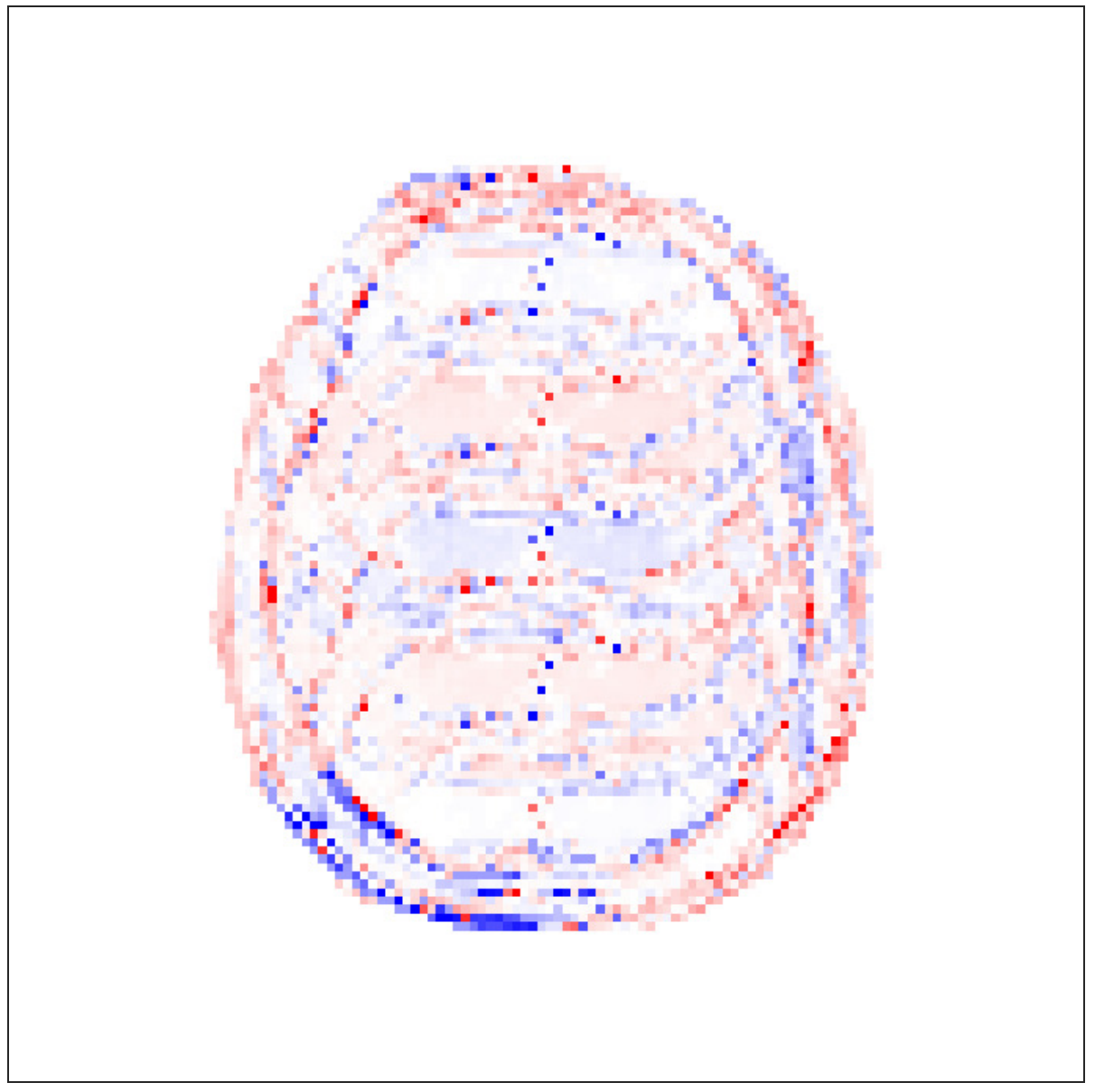}};

\draw (4.2cm, -4.3cm) node[anchor=south] {\includegraphics[width=2.1cm, height=2.1cm, clip, trim={1cm 1cm 1cm 1cm}]{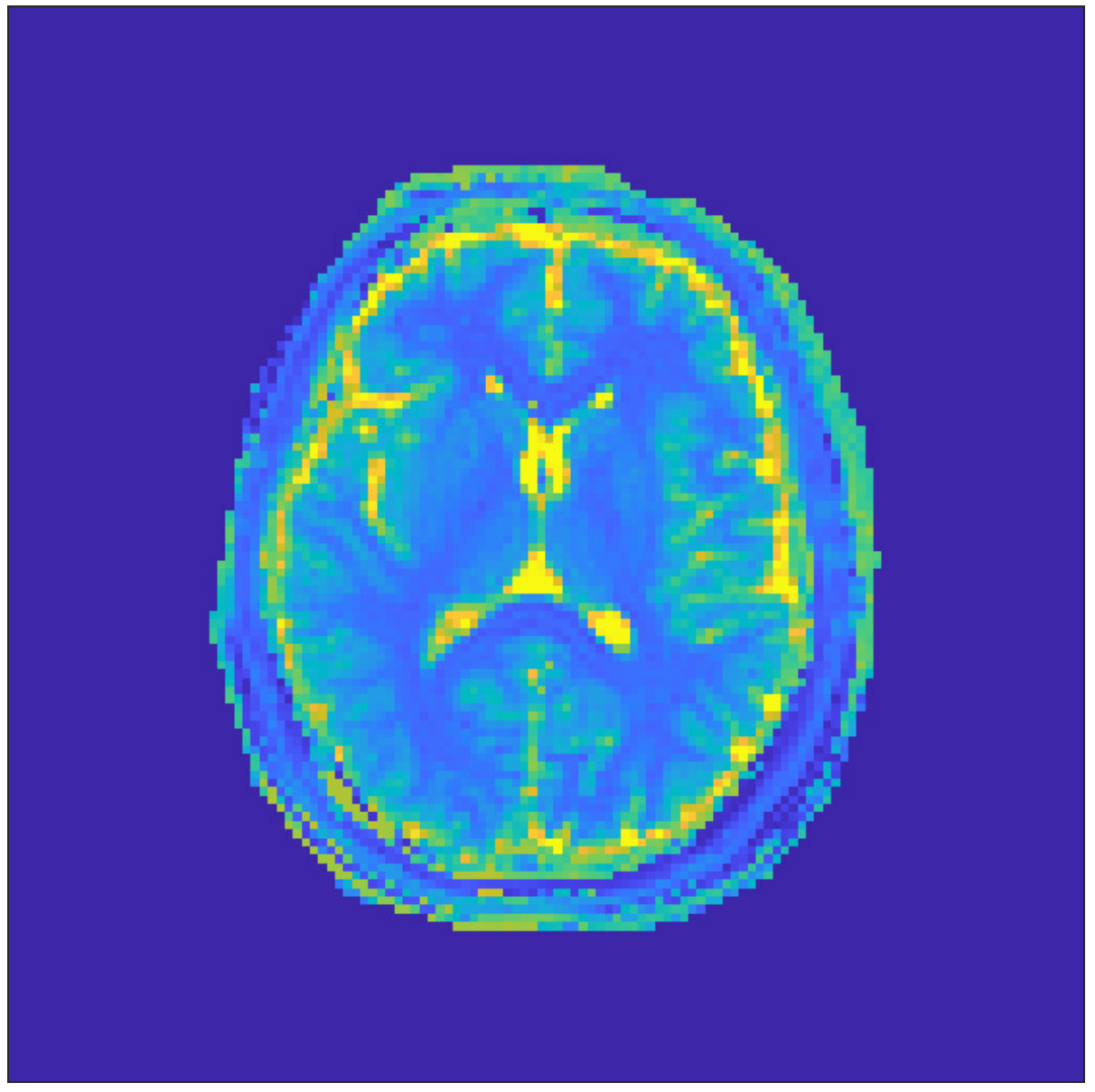}};

\draw (6.3cm, -4.3cm) node[anchor=south] {\includegraphics[width=2.1cm, height=2.1cm, clip, trim={1cm 1cm 1cm 1cm}]{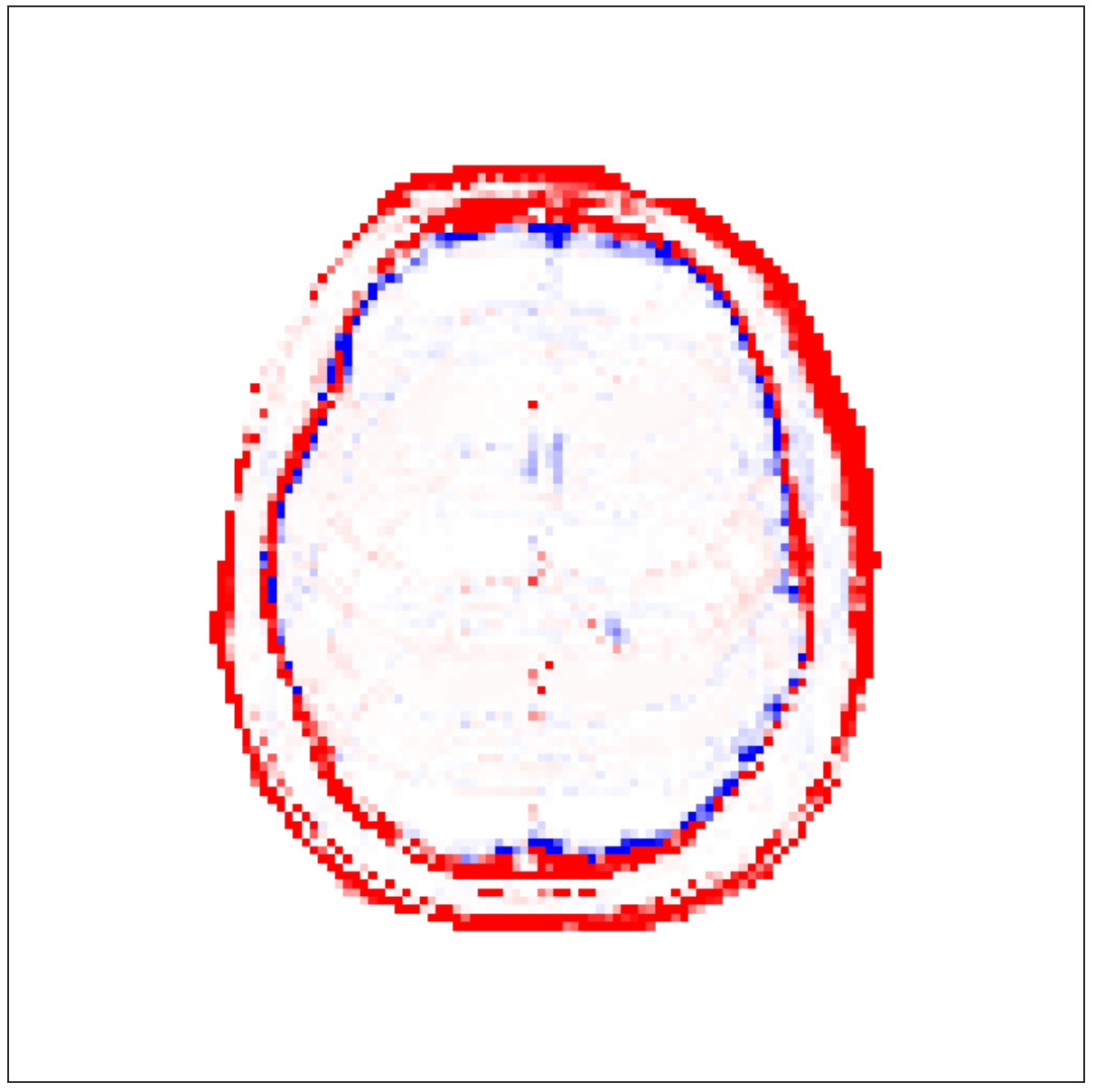}};

\draw (8.4cm, -4.3cm) node[anchor=south] {\includegraphics[width=2.1cm, height=2.1cm, clip, trim={1cm 1cm 1cm 1cm}]{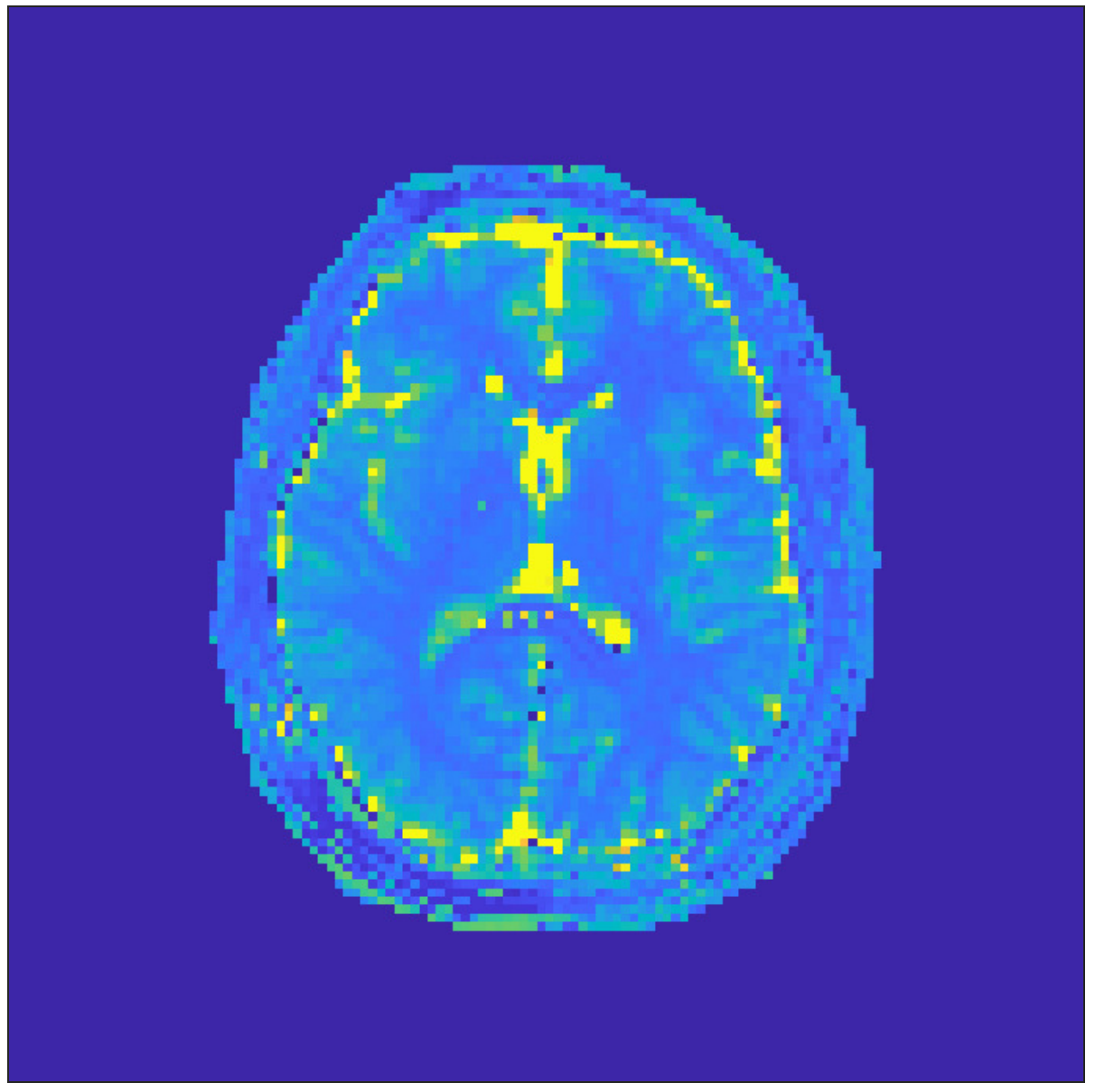}};

\draw (10.5cm, -4.3cm) node[anchor=south] {\includegraphics[width=2.1cm, height=2.1cm, clip, trim={1cm 1cm 1cm 1cm}]{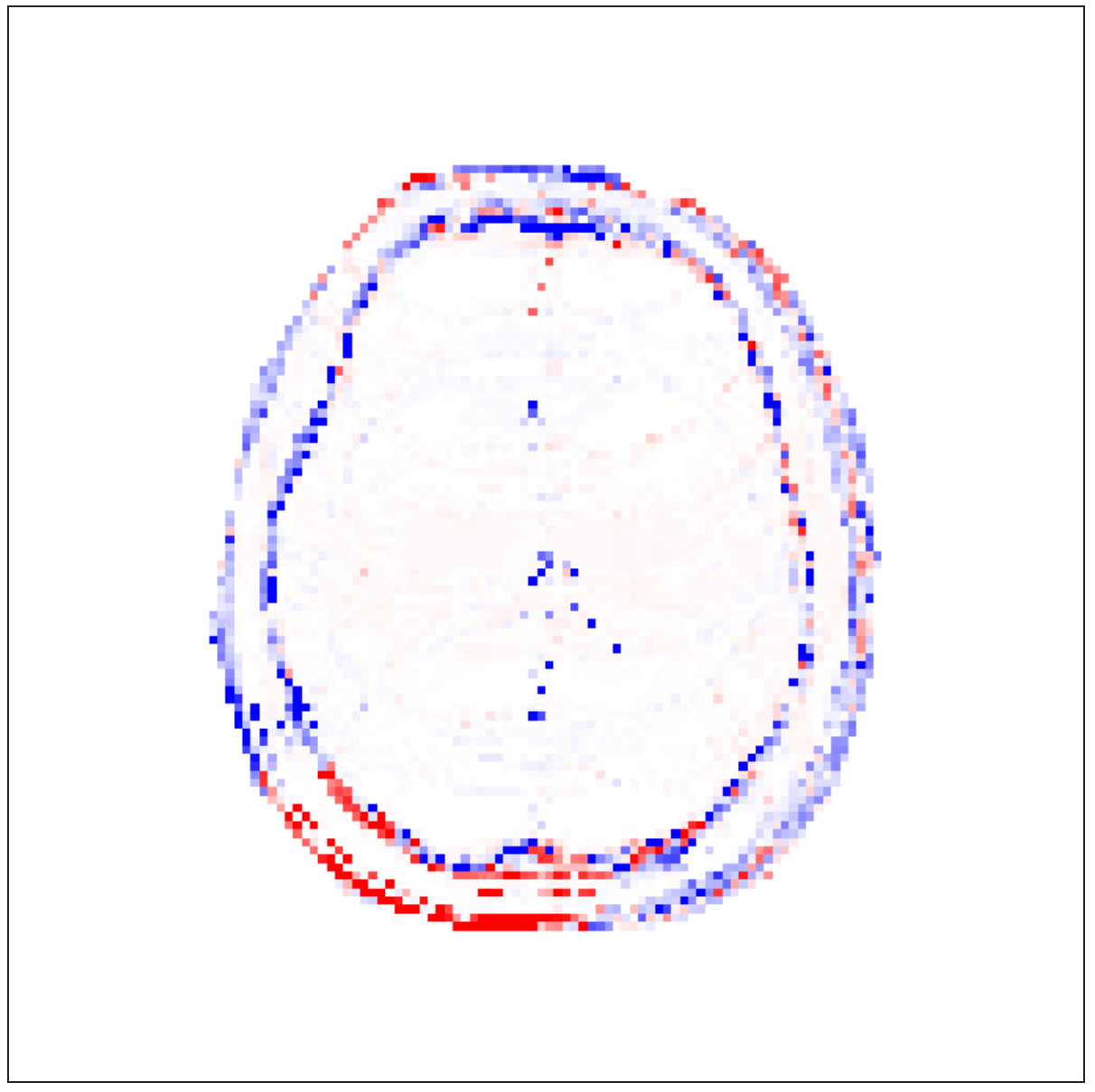}};

\draw (12.6cm, -4.3cm) node[anchor=south] {\includegraphics[width=2.1cm, height=2.1cm, clip, trim={1cm 1cm 1cm 1cm}]{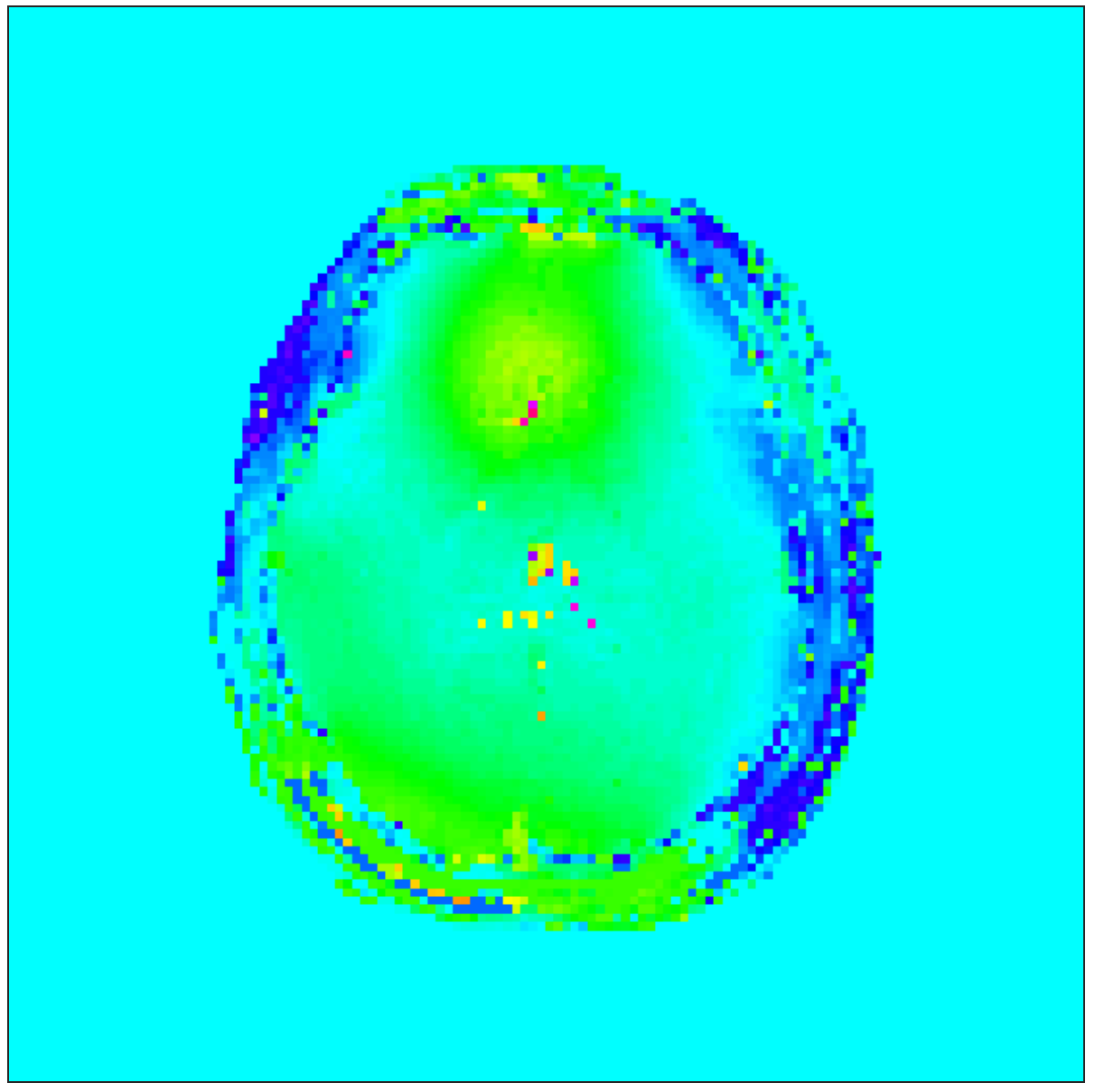}};

\draw (14.7cm, -4.3cm) node[anchor=south] {\includegraphics[width=2.1cm, height=2.1cm, clip, trim={1cm 1cm 1cm 1cm}]{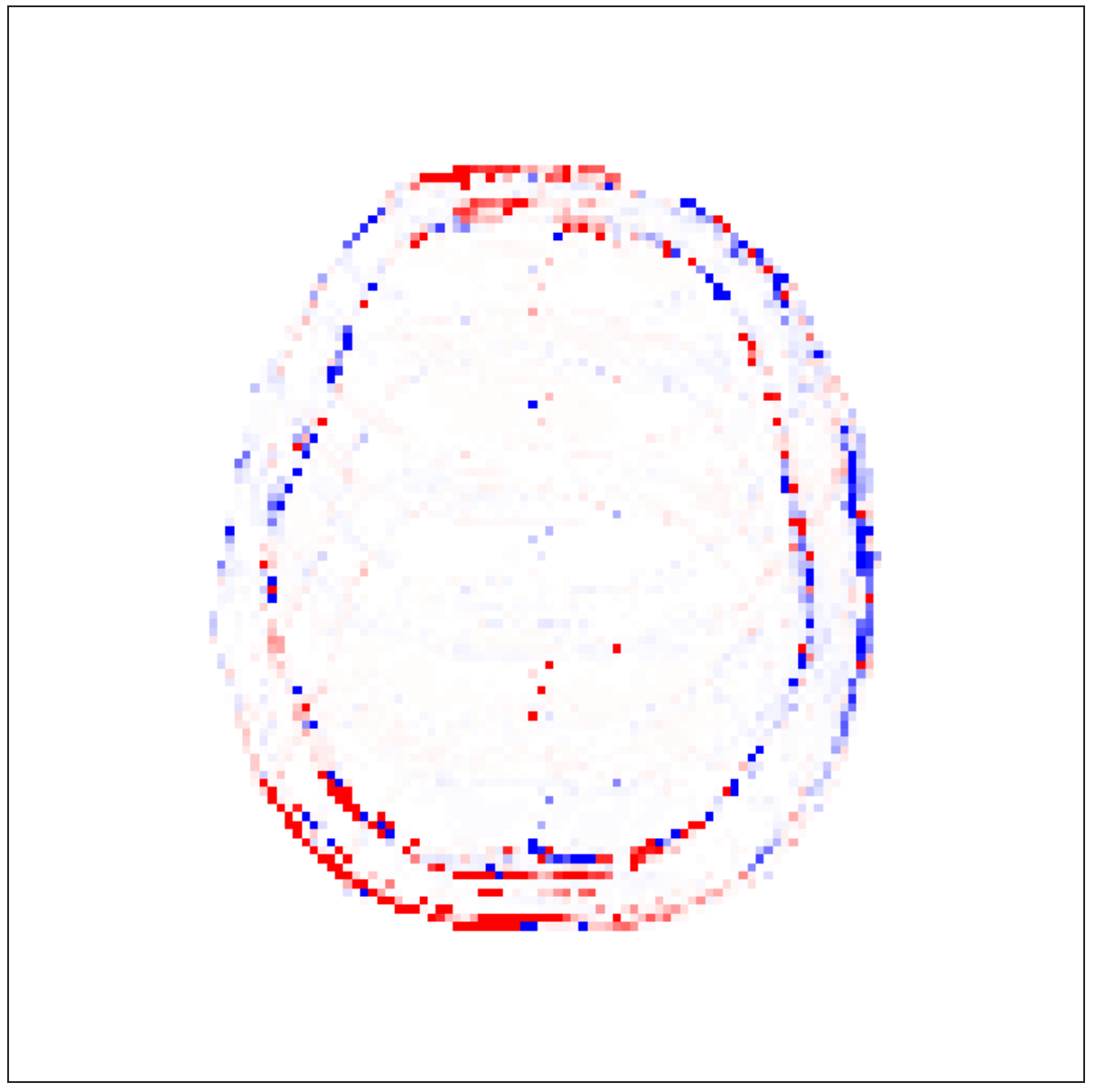}};
%%% new row
\draw (0cm, -6.45cm) node[anchor=south] {\includegraphics[width=2.1cm, height=2.1cm, clip, trim={1cm 1cm 1cm 1cm}]{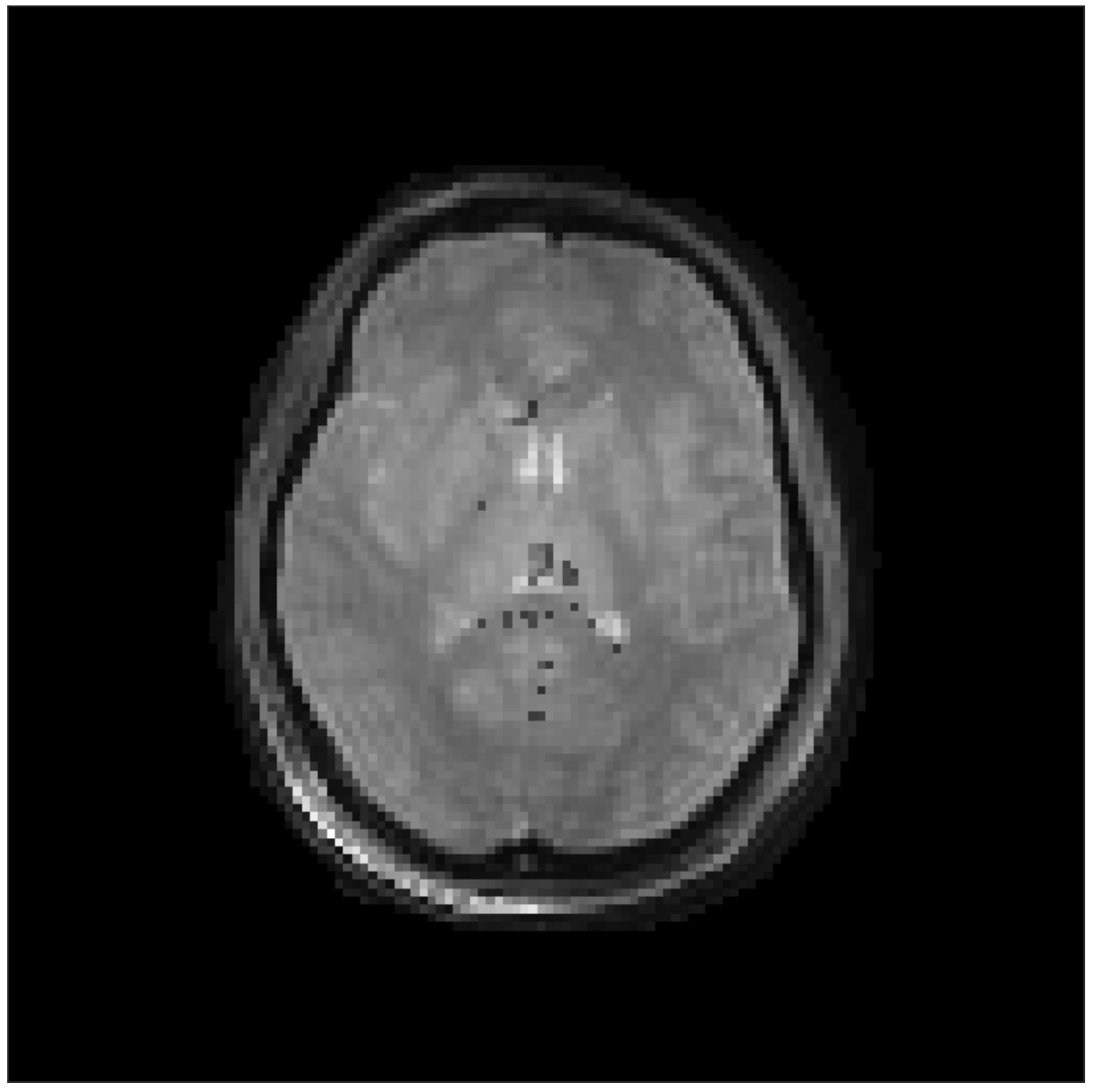}};

\draw (2.1cm, -6.45cm) node[anchor=south] {\includegraphics[width=2.1cm, height=2.1cm, clip, trim={1cm 1cm 1cm 1cm}]{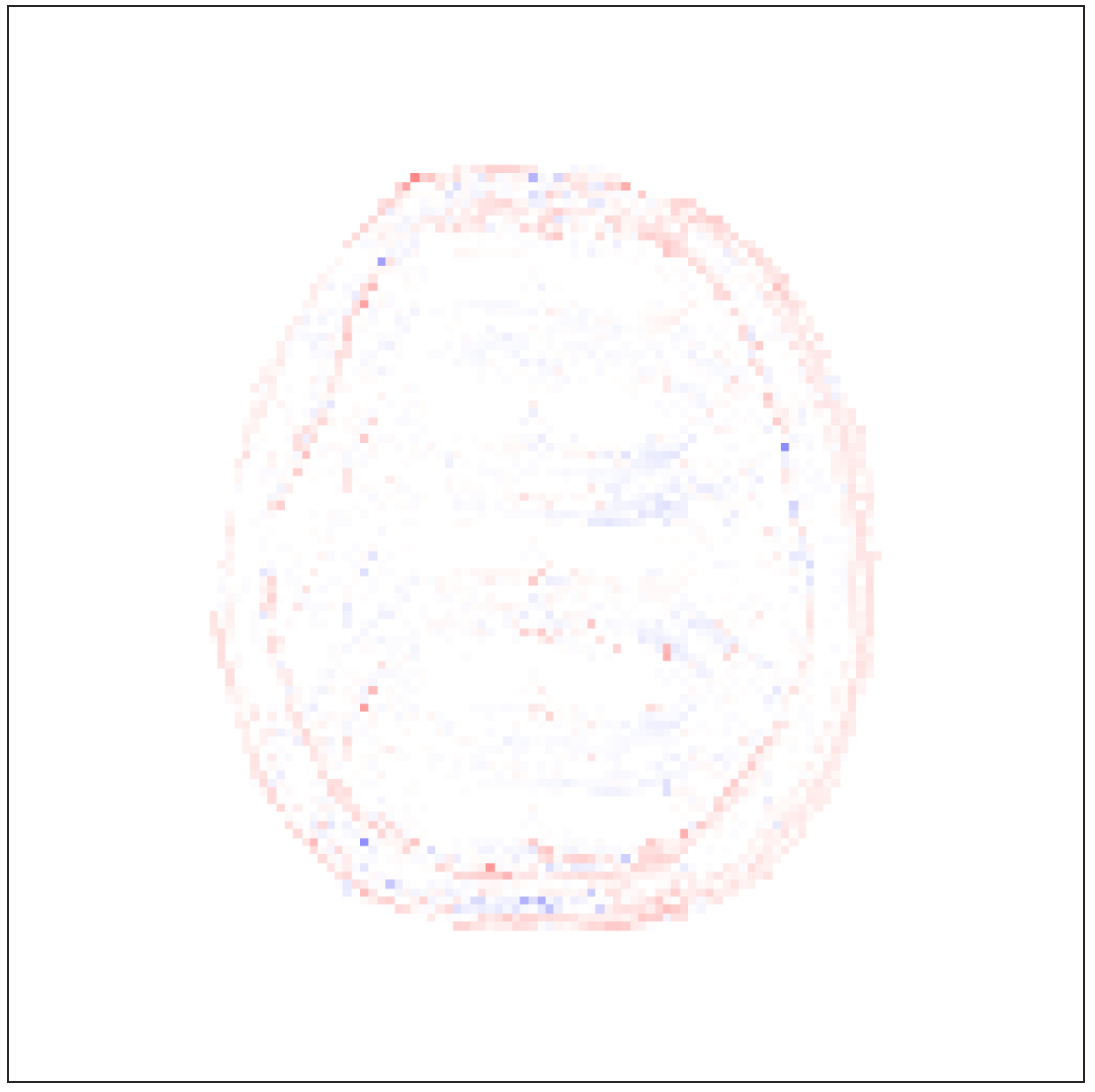}};

\draw (4.2cm, -6.45cm) node[anchor=south] {\includegraphics[width=2.1cm, height=2.1cm, clip, trim={1cm 1cm 1cm 1cm}]{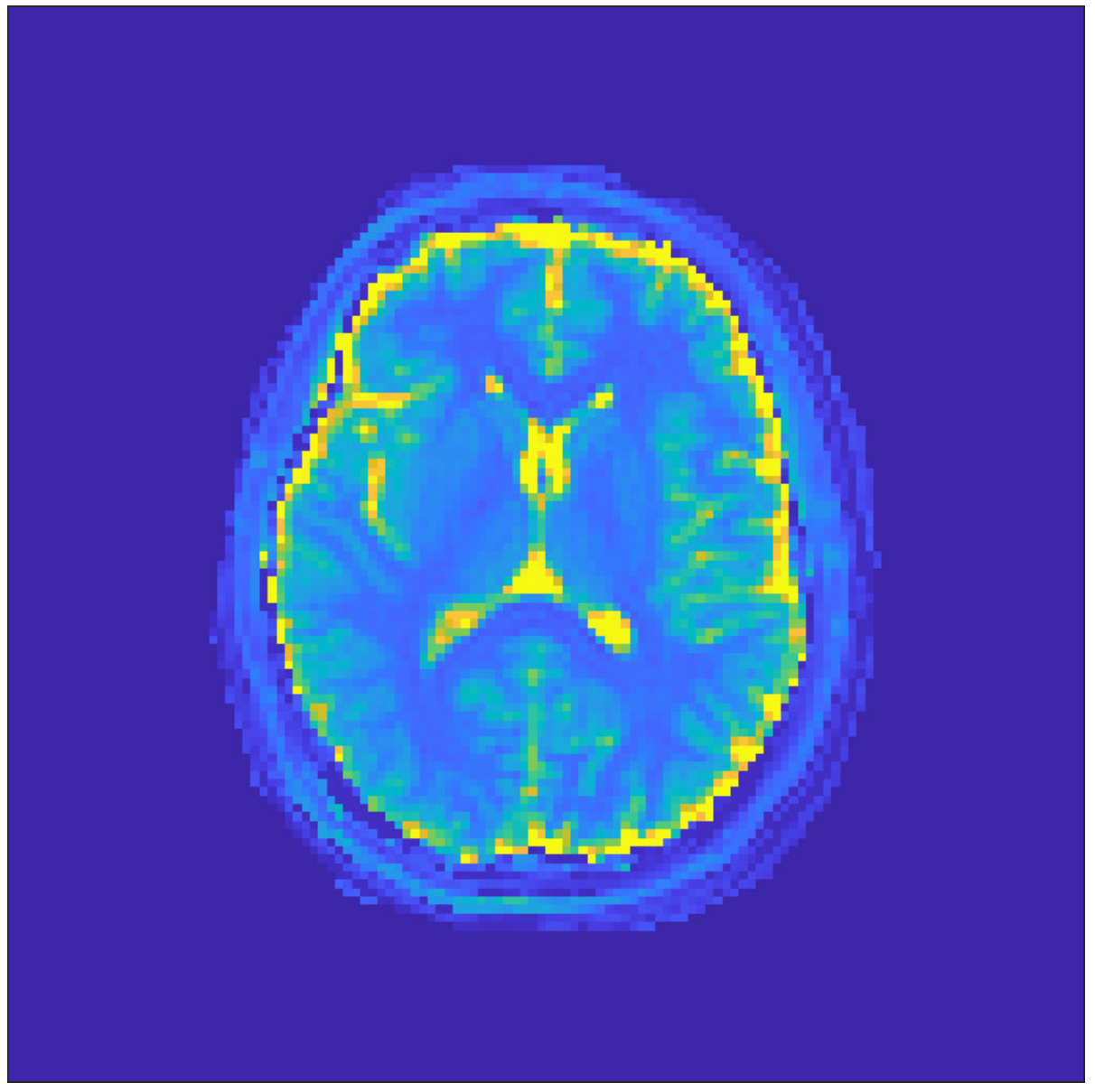}};

\draw (6.3cm, -6.45cm) node[anchor=south] {\includegraphics[width=2.1cm, height=2.1cm, clip, trim={1cm 1cm 1cm 1cm}]{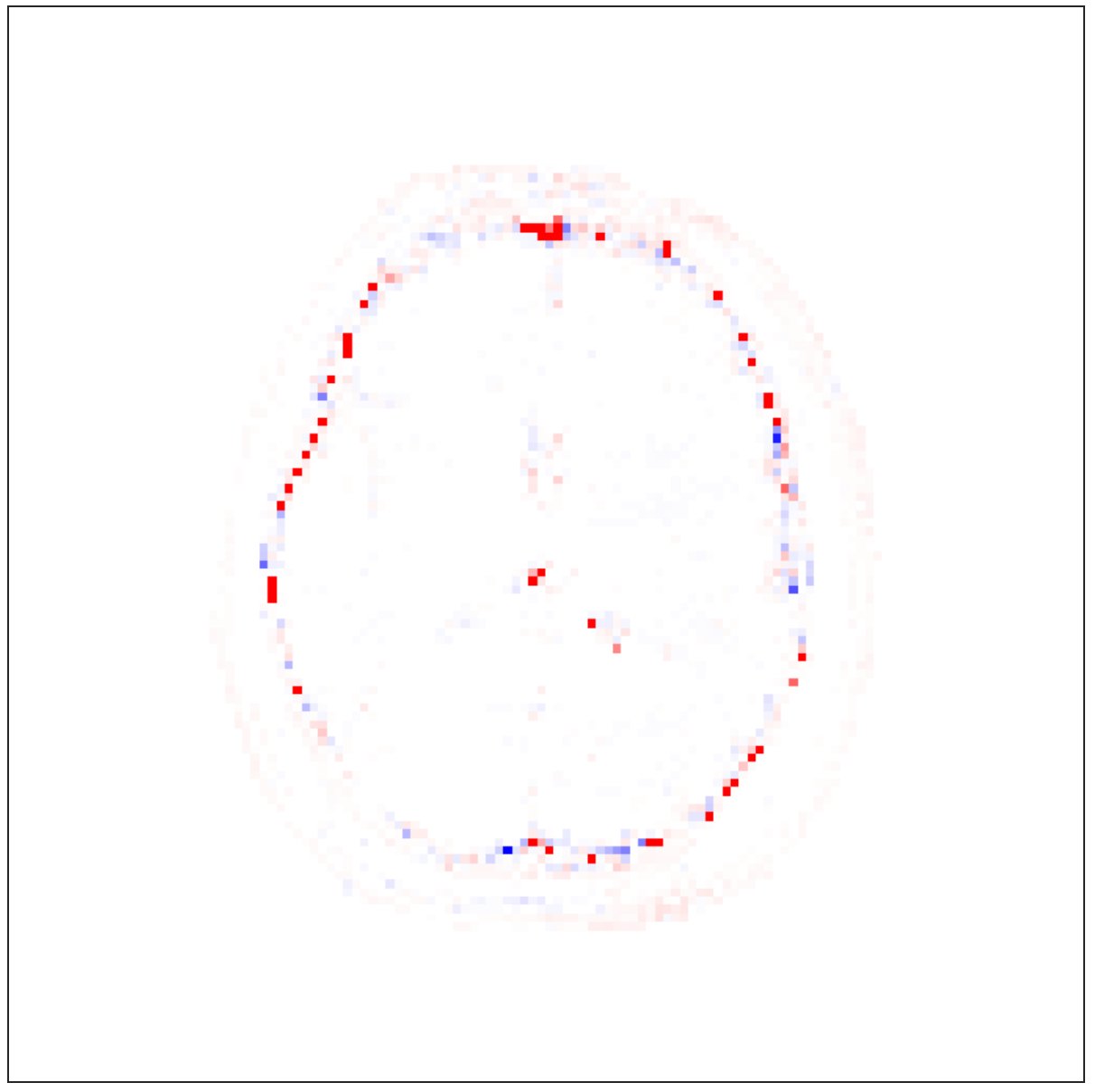}};

\draw (8.4cm, -6.45cm) node[anchor=south] {\includegraphics[width=2.1cm, height=2.1cm, clip, trim={1cm 1cm 1cm 1cm}]{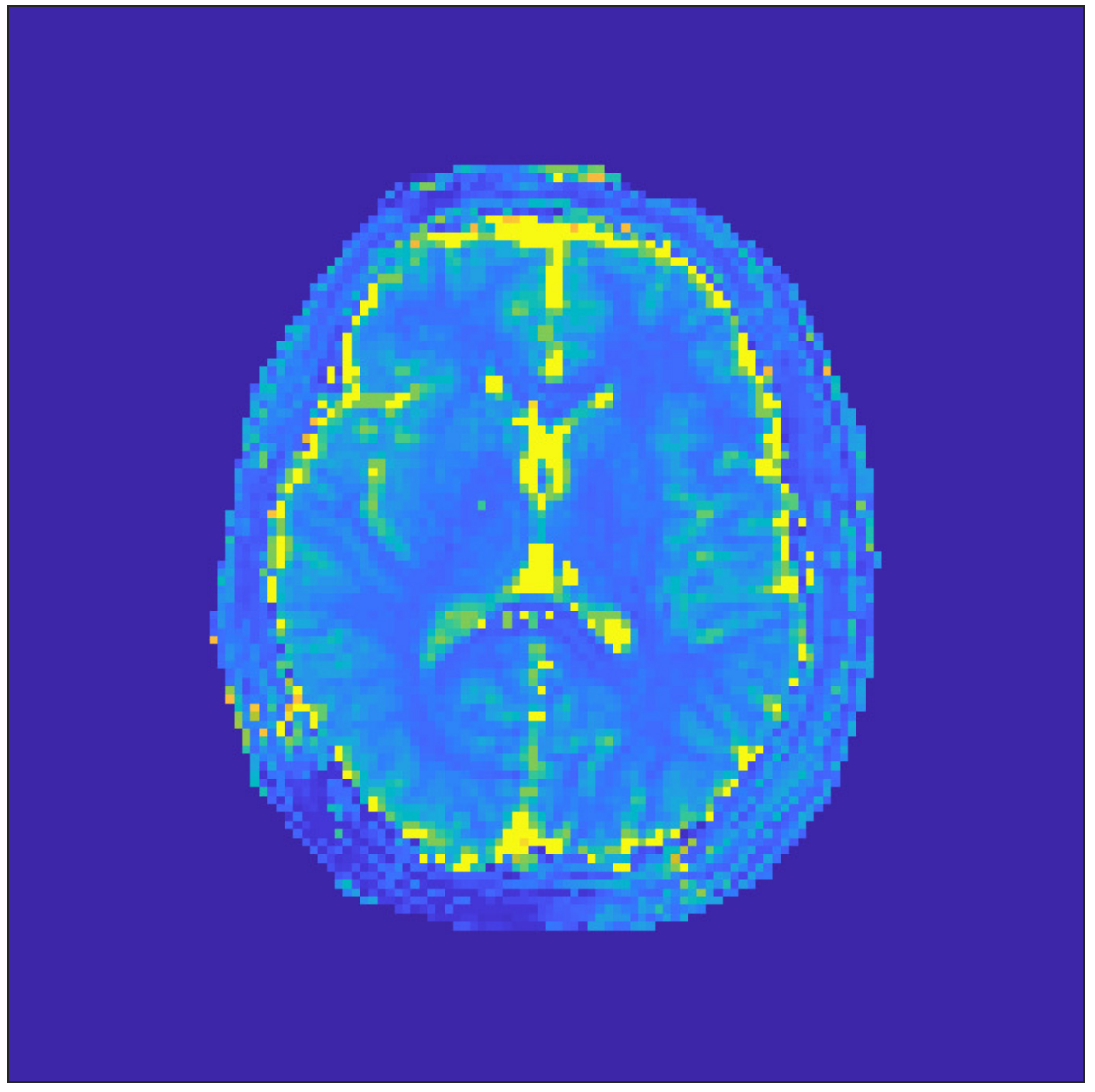}};

\draw (10.5cm, -6.45cm) node[anchor=south] {\includegraphics[width=2.1cm, height=2.1cm, clip, trim={1cm 1cm 1cm 1cm}]{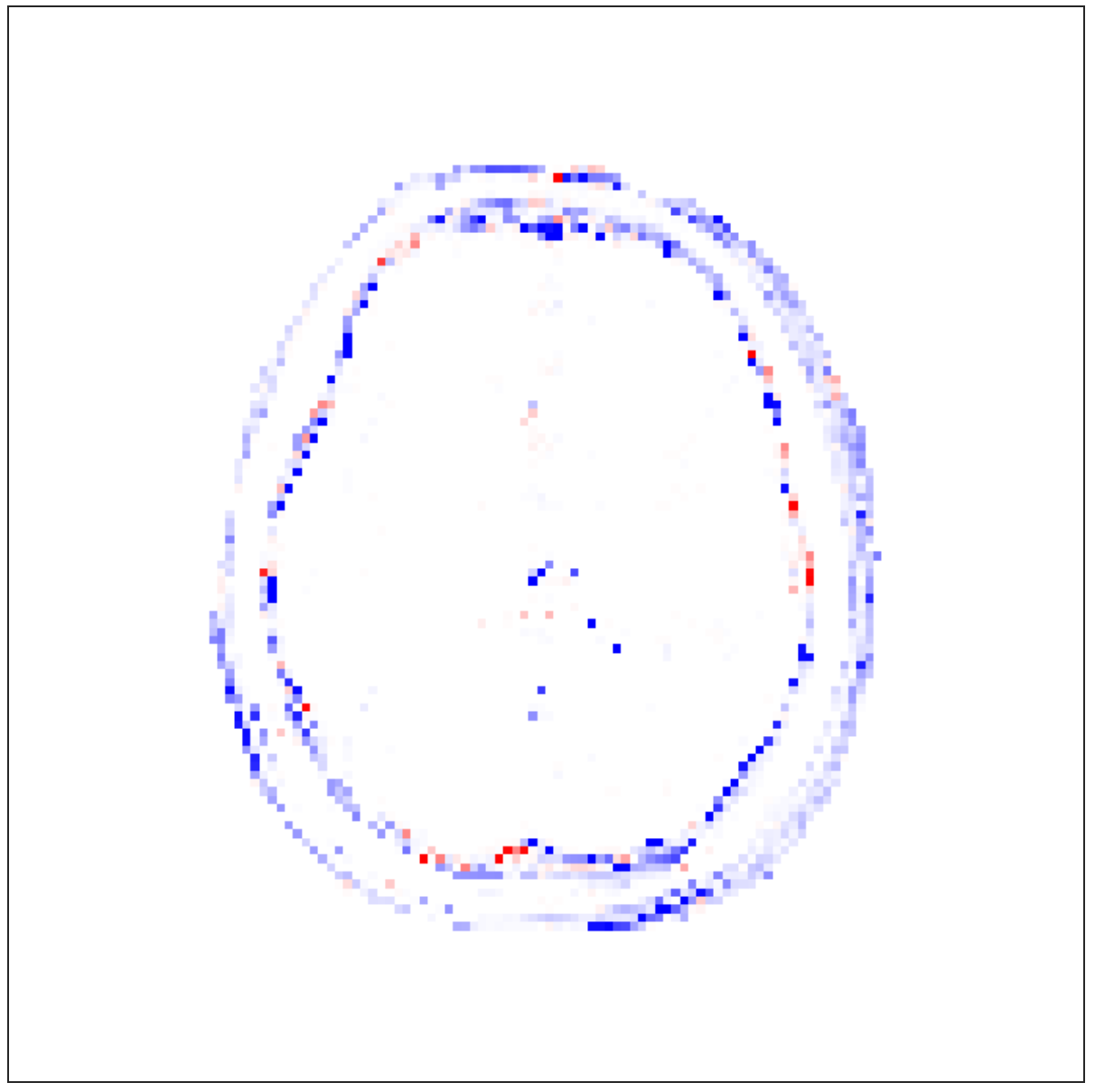}};

\draw (12.6cm, -6.45cm) node[anchor=south] {\includegraphics[width=2.1cm, height=2.1cm, clip, trim={1cm 1cm 1cm 1cm}]{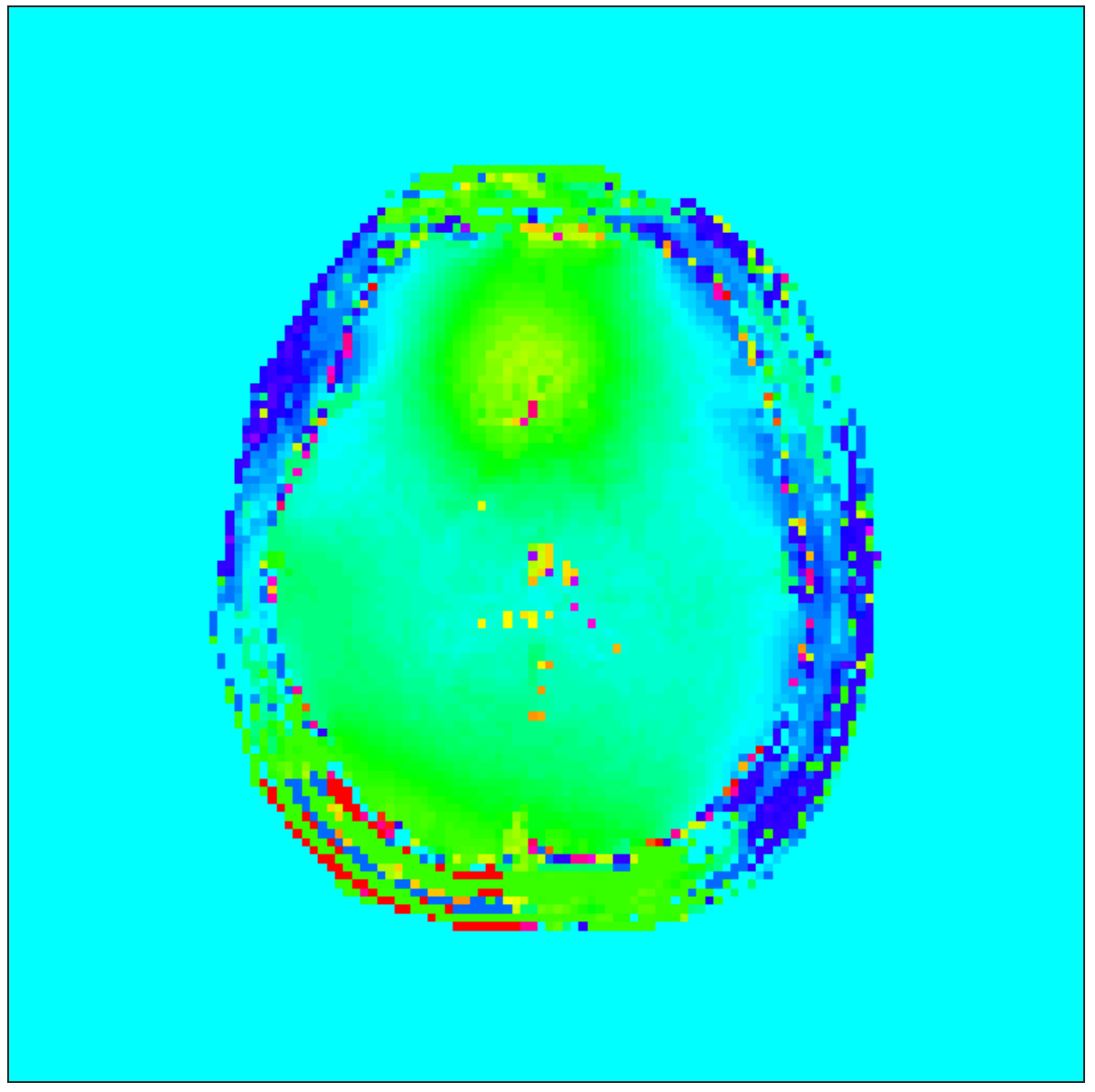}};

\draw (14.7cm, -6.45cm) node[anchor=south] {\includegraphics[width=2.1cm, height=2.1cm, clip, trim={1cm 1cm 1cm 1cm}]{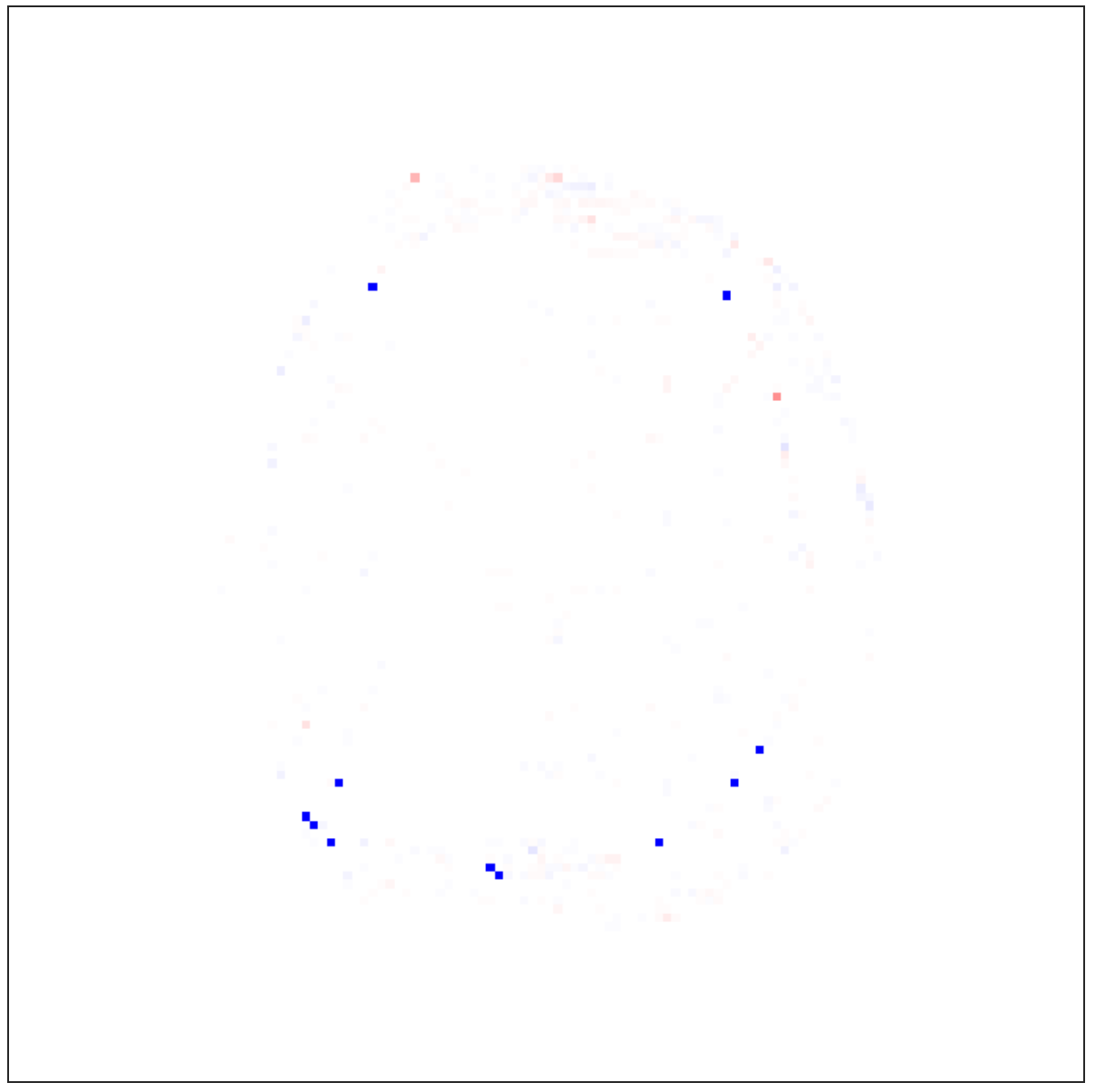}};
%%% new row
\draw (0cm, -8.6cm) node[anchor=south] {\includegraphics[width=2.1cm, height=2.1cm, clip, trim={1cm 1cm 1cm 1cm}]{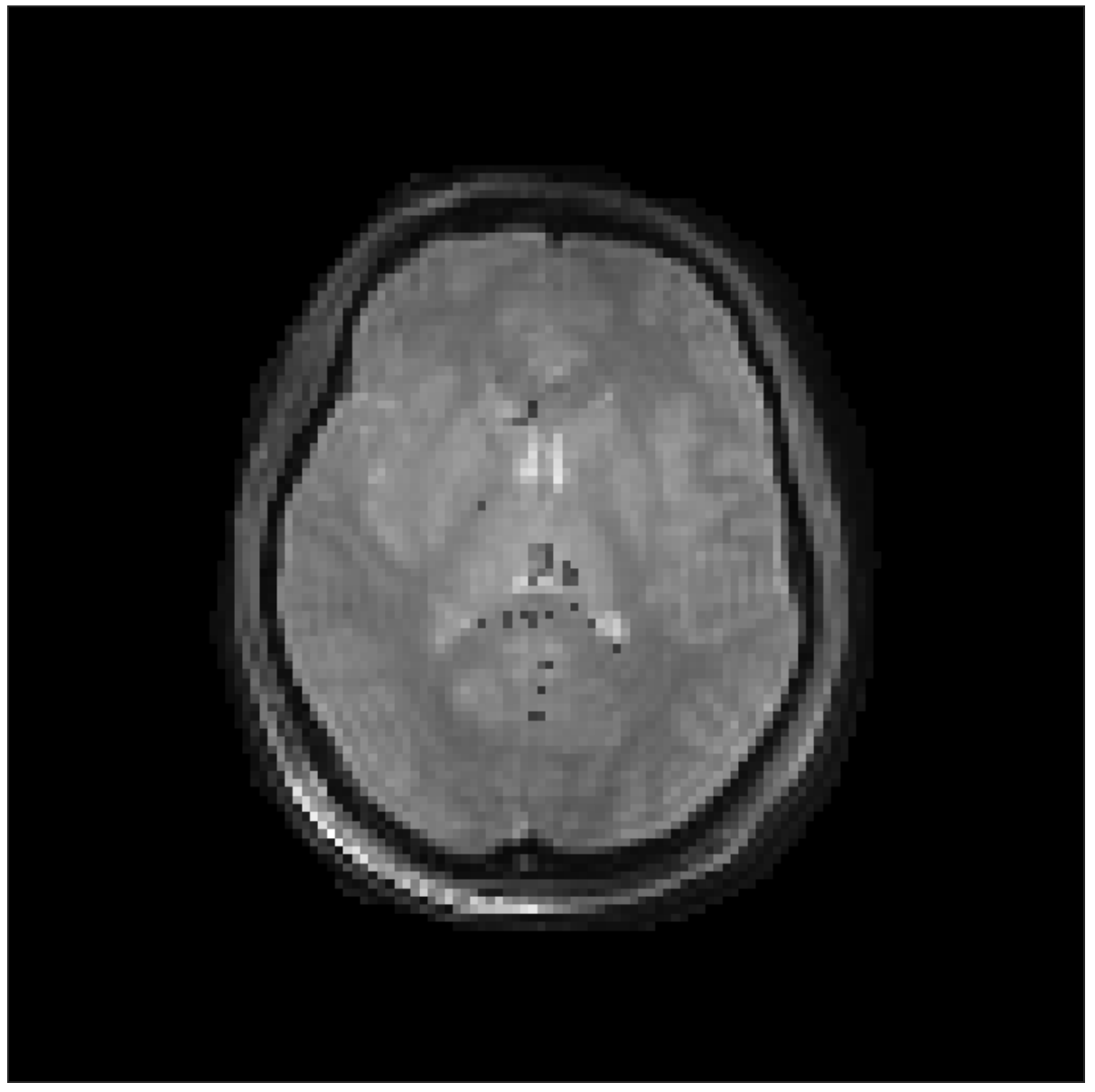}};

\draw (2.1cm, -8.6cm) node[anchor=south] {\includegraphics[width=2.1cm, height=2.1cm, clip, trim={1cm 1cm 1cm 1cm}]{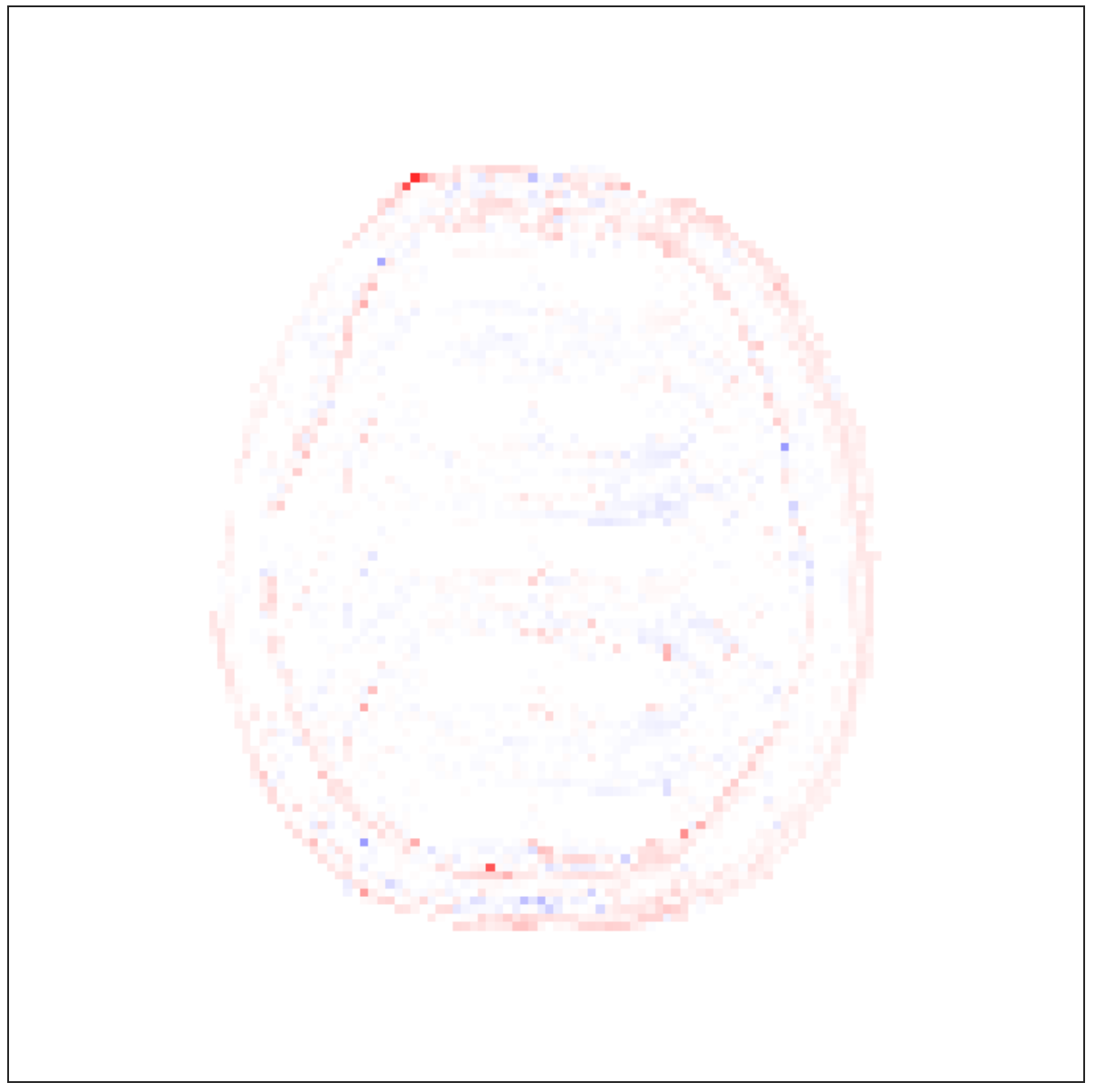}};

\draw (4.2cm, -8.6cm) node[anchor=south] {\includegraphics[width=2.1cm, height=2.1cm, clip, trim={1cm 1cm 1cm 1cm}]{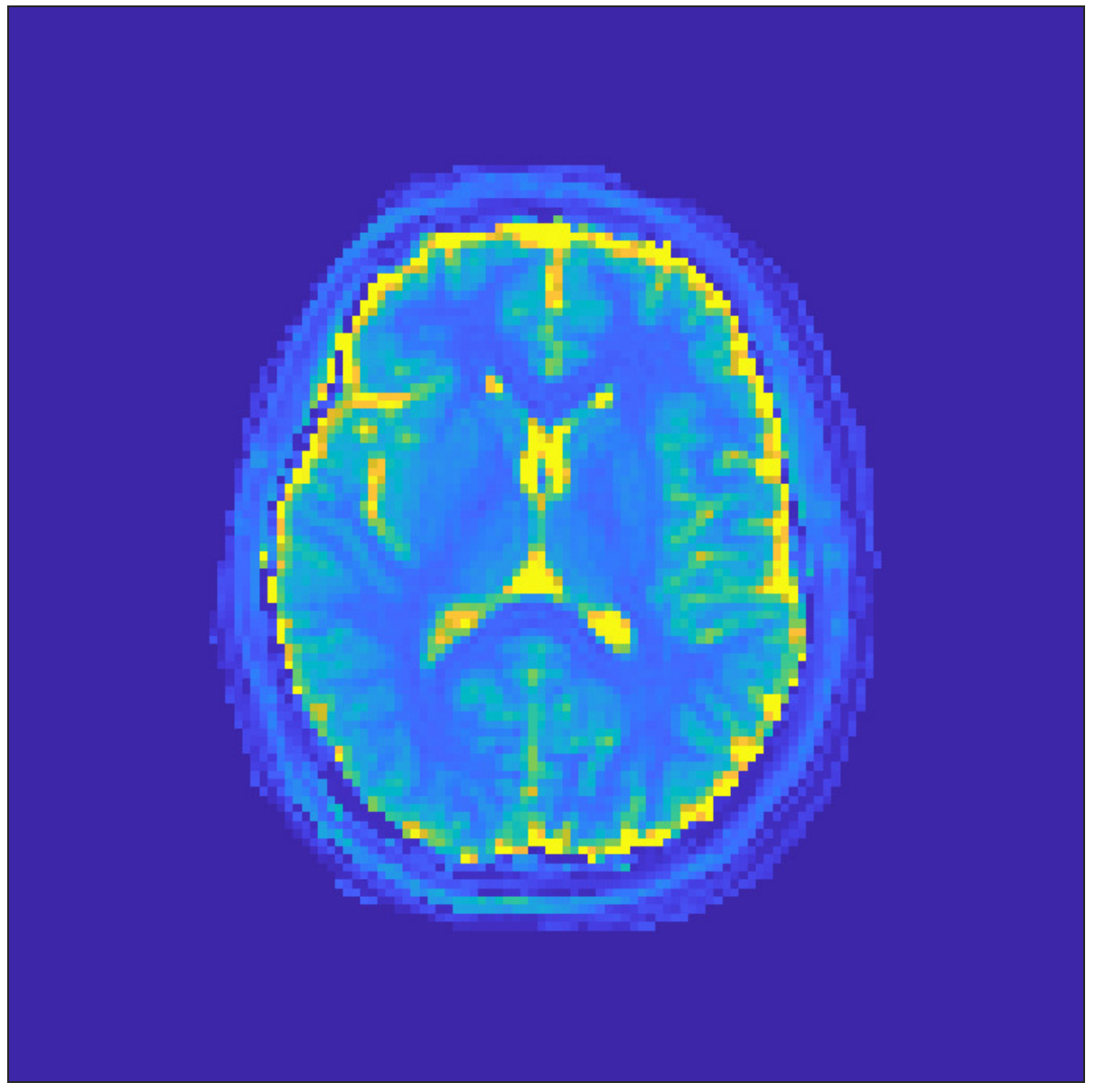}};

\draw (6.3cm, -8.6cm) node[anchor=south] {\includegraphics[width=2.1cm, height=2.1cm, clip, trim={1cm 1cm 1cm 1cm}]{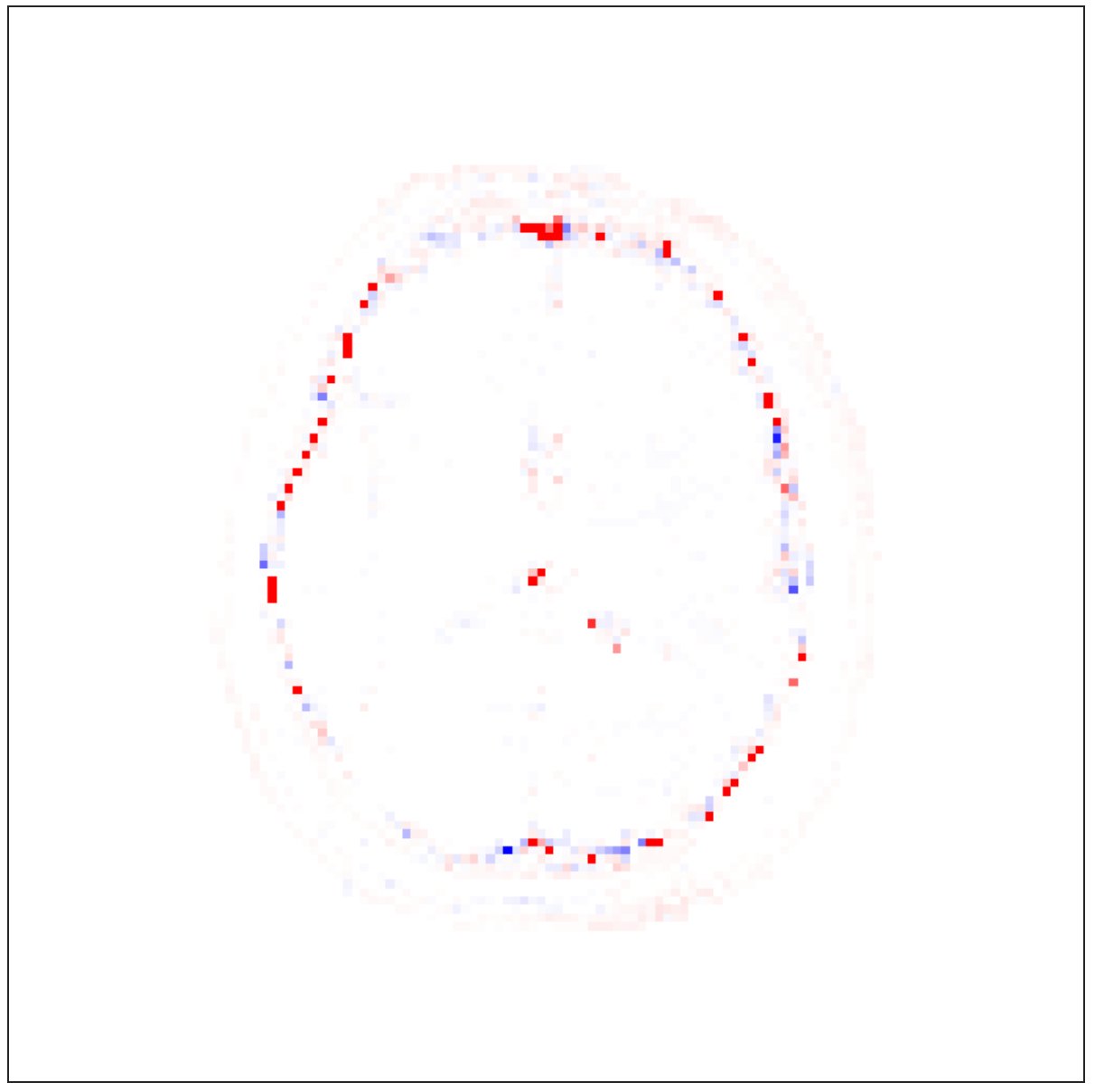}};

\draw (8.4cm, -8.6cm) node[anchor=south] {\includegraphics[width=2.1cm, height=2.1cm, clip, trim={1cm 1cm 1cm 1cm}]{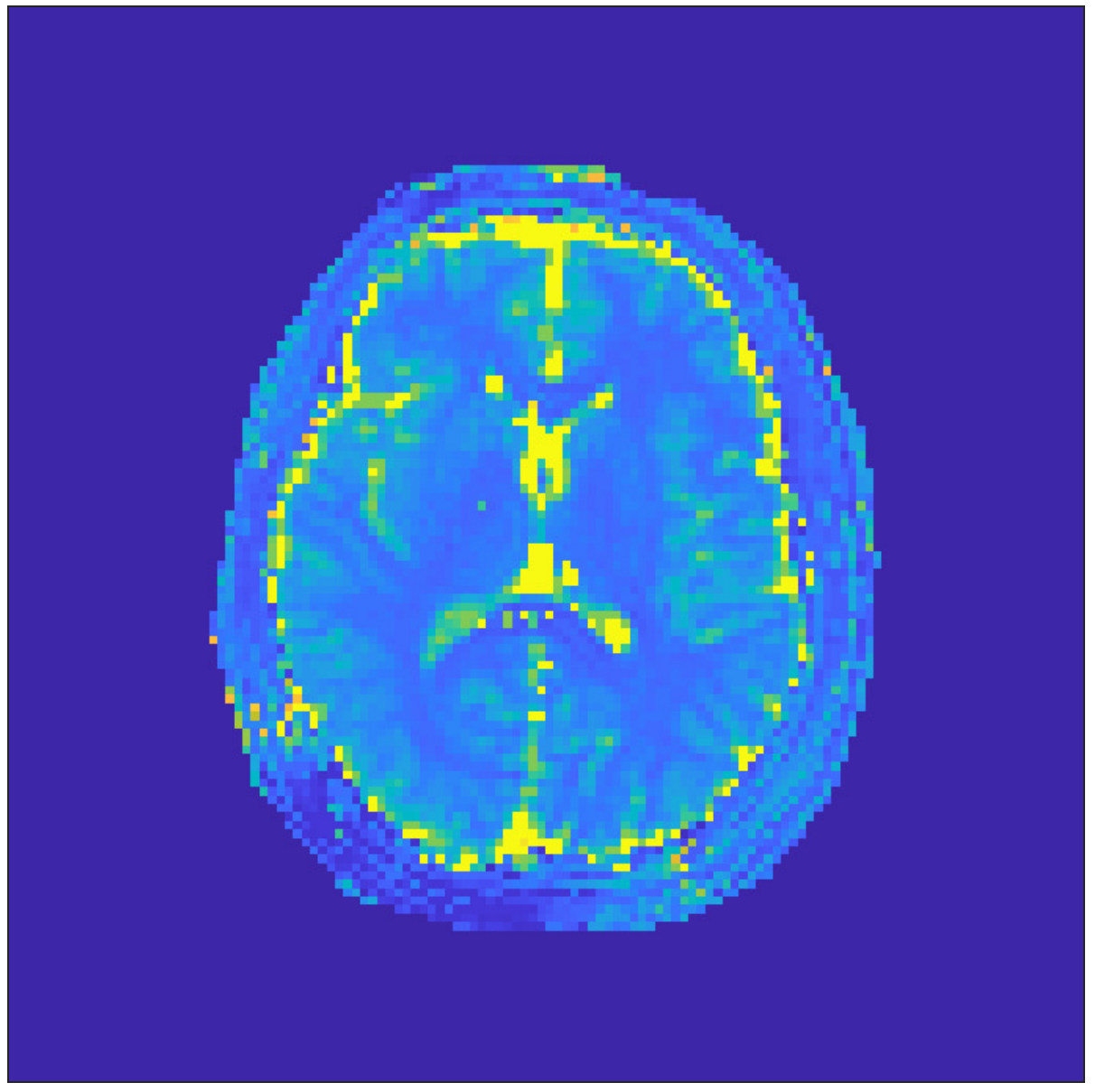}};

\draw (10.5cm, -8.6cm) node[anchor=south] {\includegraphics[width=2.1cm, height=2.1cm, clip, trim={1cm 1cm 1cm 1cm}]{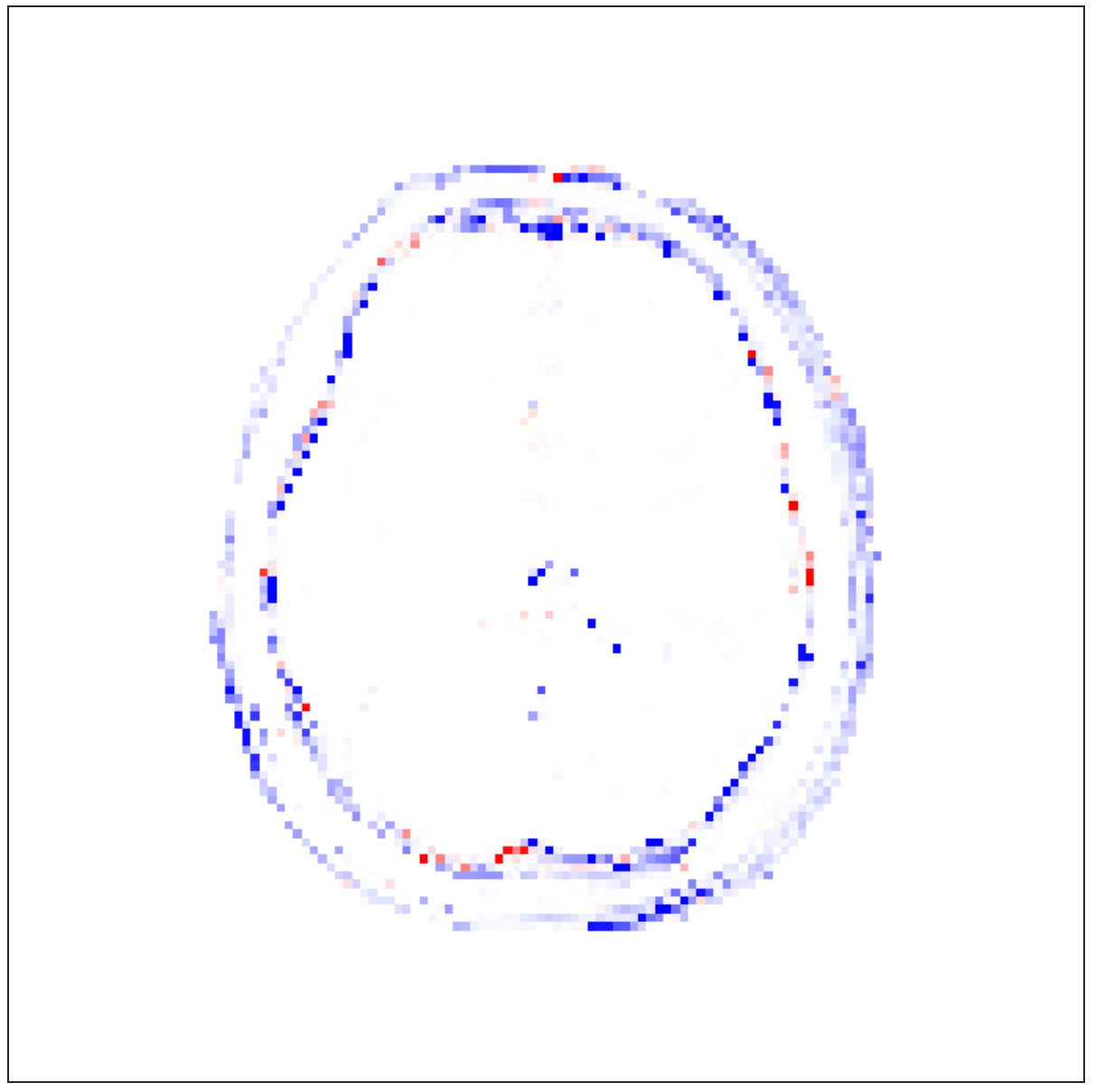}};

\draw (12.6cm, -8.6cm) node[anchor=south] {\includegraphics[width=2.1cm, height=2.1cm, clip, trim={1cm 1cm 1cm 1cm}]{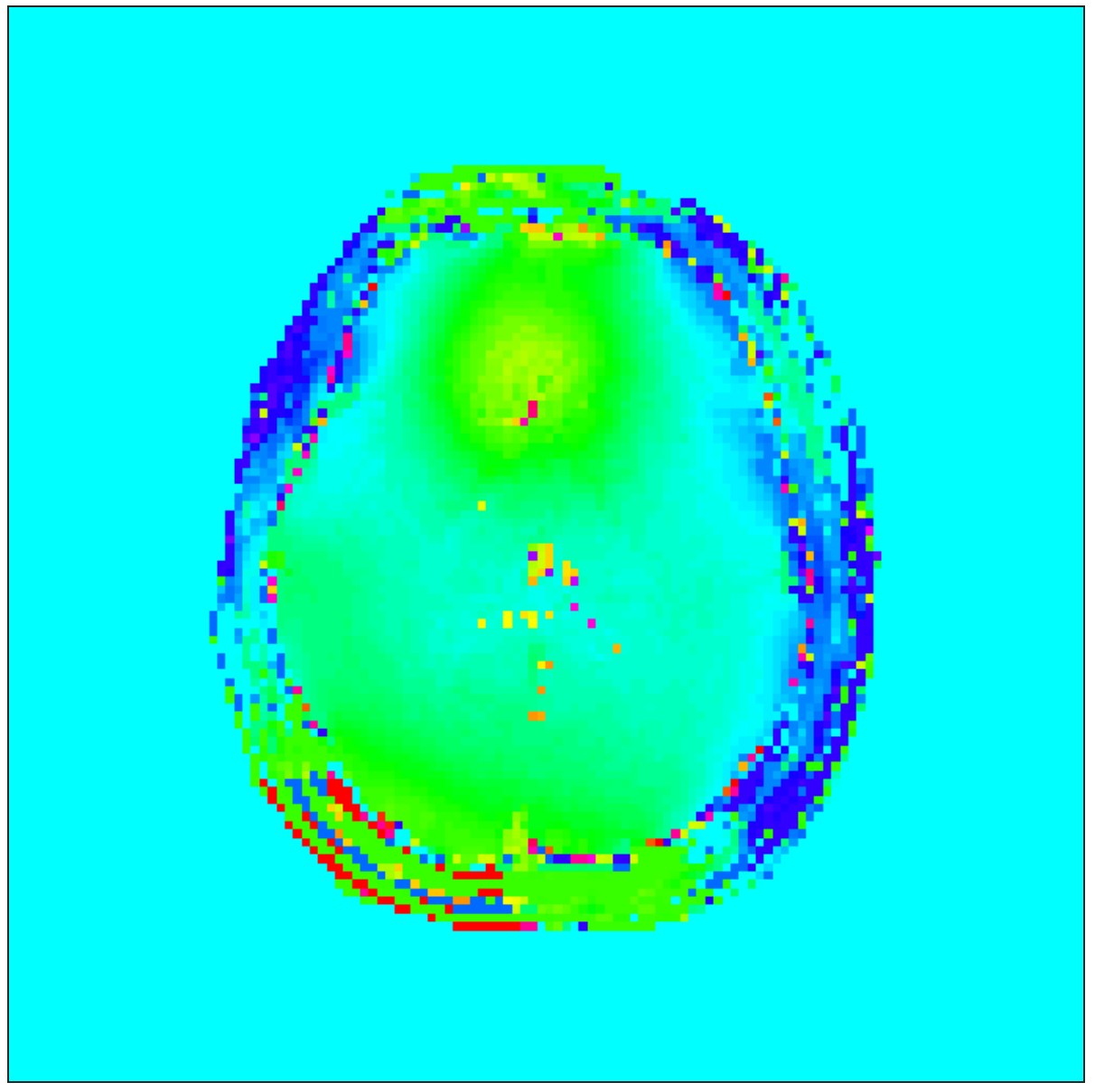}};

\draw (14.7cm, -8.6cm) node[anchor=south] {\includegraphics[width=2.1cm, height=2.1cm, clip, trim={1cm 1cm 1cm 1cm}]{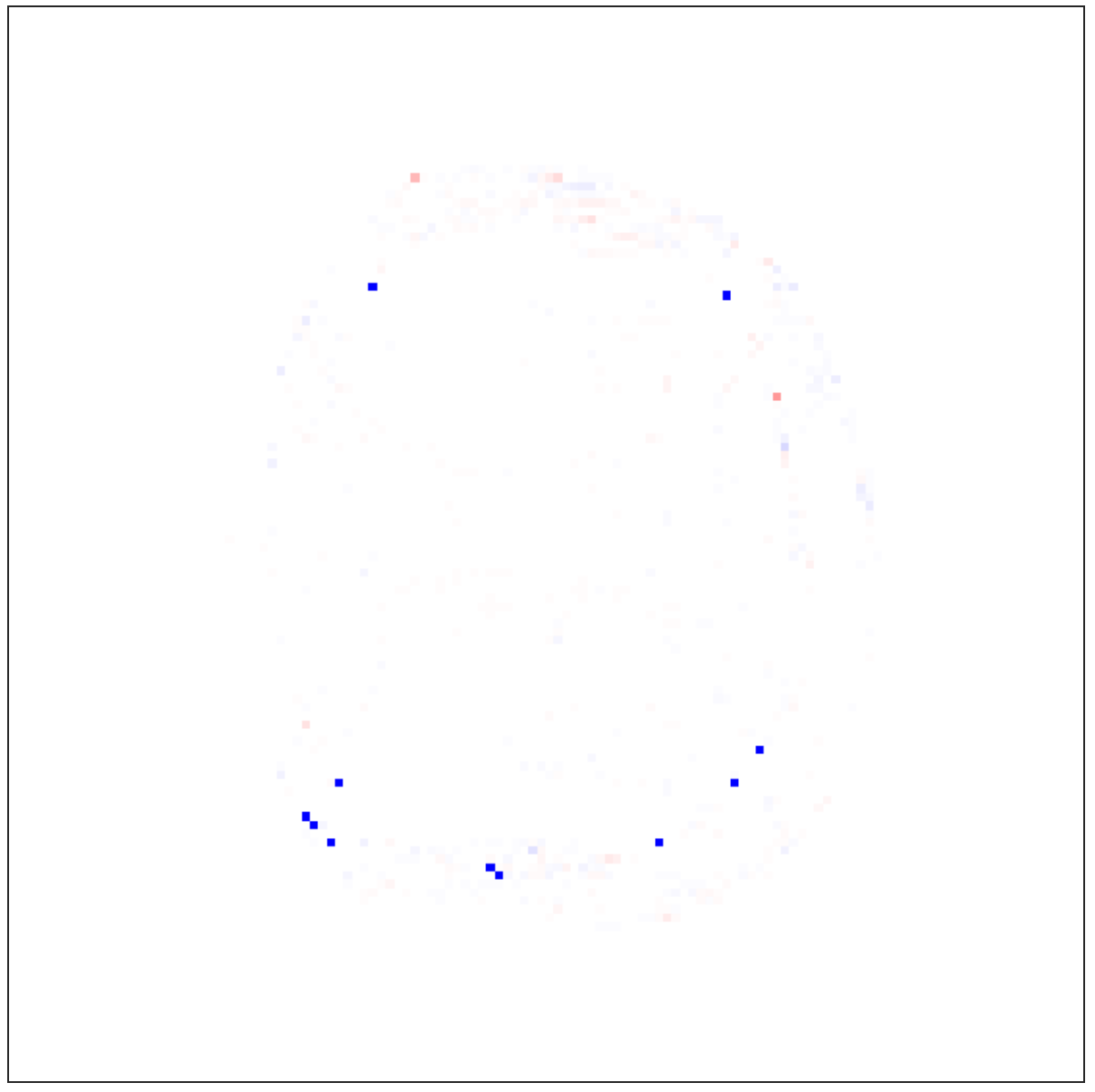}};
%%% colorbars
\draw (0cm, -8.85cm) node[anchor=south] {\includegraphics[width=2cm, clip, trim={0cm 0cm 0cm 0cm}]{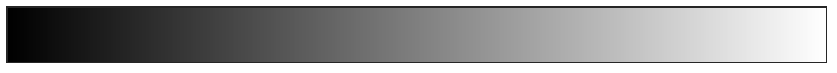}};

\draw (2.1cm, -8.85cm) node[anchor=south] {\includegraphics[width=2cm, clip, trim={0cm 0cm 0cm 0cm}]{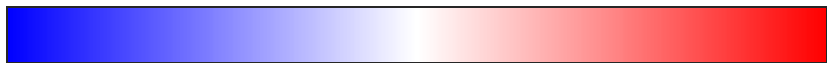}};

\draw (4.2cm, -8.85cm) node[anchor=south] {\includegraphics[width=2cm, clip, trim={0cm 0cm 0cm 0cm}]{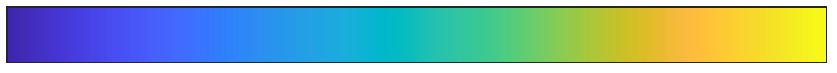}};

\draw (6.3cm, -8.85cm) node[anchor=south] {\includegraphics[width=2cm, clip, trim={0cm 0cm 0cm 0cm}]{figures/colorbars/colorbar_redblue.pdf}};

\draw (8.4cm, -8.85cm) node[anchor=south] {\includegraphics[width=2cm, clip, trim={0cm 0cm 0cm 0cm}]{figures/colorbars/colorbar_parula.pdf}};

\draw (10.5cm, -8.85cm) node[anchor=south] {\includegraphics[width=2cm, clip, trim={0cm 0cm 0cm 0cm}]{figures/colorbars/colorbar_redblue.pdf}};

\draw (12.6cm, -8.85cm) node[anchor=south] {\includegraphics[width=2cm, clip, trim={0cm 0cm 0cm 0cm}]{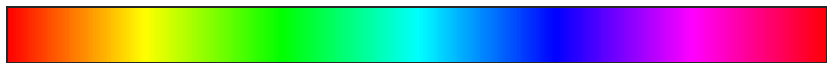}};

\draw (14.7cm, -8.85cm) node[anchor=south] {\includegraphics[width=2cm, clip, trim={0cm 0cm 0cm 0cm}]{figures/colorbars/colorbar_redblue.pdf}};
%%% colorbar labels
\draw (-0.92cm, -9.1cm) node[anchor=south] {\begin{scriptsize} 0 \end{scriptsize}};
\draw (-0.02cm, -9.1cm) node[anchor=south] {\begin{scriptsize} 0.5 \end{scriptsize}};
\draw (0.9cm, -9.1cm) node[anchor=south] {\begin{scriptsize} 1 \end{scriptsize}};

\draw (1.28cm, -9.1cm) node[anchor=south] {\begin{scriptsize} -0.1 \end{scriptsize}};
\draw (2.08cm, -9.1cm) node[anchor=south] {\begin{scriptsize} 0 \end{scriptsize}};
\draw (2.9cm, -9.1cm) node[anchor=south] {\begin{scriptsize} 0.1 \end{scriptsize}};

\draw (3.28cm, -9.1cm) node[anchor=south] {\begin{scriptsize} 0 \end{scriptsize}};
\draw (4.18cm, -9.1cm) node[anchor=south] {\begin{scriptsize} 1.5 \end{scriptsize}};
\draw (5.1cm, -9.1cm) node[anchor=south] {\begin{scriptsize} 3 \end{scriptsize}};

\draw (5.38cm, -9.1cm) node[anchor=south] {\begin{scriptsize} -1 \end{scriptsize}};
\draw (6.28cm, -9.1cm) node[anchor=south] {\begin{scriptsize} 0 \end{scriptsize}};
\draw (7.2cm, -9.1cm) node[anchor=south] {\begin{scriptsize} 1 \end{scriptsize}};

\draw (7.48cm, -9.1cm) node[anchor=south] {\begin{scriptsize} 0 \end{scriptsize}};
\draw (8.36cm, -9.1cm) node[anchor=south] {\begin{scriptsize} 0.15 \end{scriptsize}};
\draw (9.2cm, -9.1cm) node[anchor=south] {\begin{scriptsize} 0.3 \end{scriptsize}};

\draw (9.7cm, -9.1cm) node[anchor=south] {\begin{scriptsize} -0.1 \end{scriptsize}};
\draw (10.48cm, -9.1cm) node[anchor=south] {\begin{scriptsize} 0 \end{scriptsize}};
\draw (11.3cm, -9.1cm) node[anchor=south] {\begin{scriptsize} 0.1 \end{scriptsize}};

\draw (11.77cm, -9.1cm) node[anchor=south] {\begin{scriptsize} -50 \end{scriptsize}};
\draw (12.58cm, -9.1cm) node[anchor=south] {\begin{scriptsize} 0 \end{scriptsize}};
\draw (13.43cm, -9.1cm) node[anchor=south] {\begin{scriptsize} 50 \end{scriptsize}};

\draw (13.86cm, -9.1cm) node[anchor=south] {\begin{scriptsize} -10 \end{scriptsize}};
\draw (14.7cm, -9.1cm) node[anchor=south] {\begin{scriptsize} 0 \end{scriptsize}};
\draw (15.51cm, -9.1cm) node[anchor=south] {\begin{scriptsize} 10 \end{scriptsize}};
%%% column labels
\draw (0cm, -9.75cm) node[anchor=south] {$\rho$ (a.u.)};
\draw (2.1cm, -9.75cm) node[anchor=south] {$\rho$ error};
\draw (4.2cm, -9.75cm) node[anchor=south] {$T_1$ (s)};
\draw (6.3cm, -9.75cm) node[anchor=south] {$T_1$ error};
\draw (8.4cm, -9.75cm) node[anchor=south] {$T_2$ (s)};
\draw (10.5cm, -9.75cm) node[anchor=south] {$T_2$ error};
\draw (12.6cm, -9.75cm) node[anchor=south] {$\omega$ (Hz)};
\draw (14.7cm, -9.75cm) node[anchor=south] {$\omega$ error};
%%% row labels
\draw (-1cm, 1.15cm) node[anchor=south,rotate=90] {BLIP};
\draw (-1cm, -1cm) node[anchor=south,rotate=90] {FINE};
\draw (-1cm, -3.15cm) node[anchor=south,rotate=90] {C2F};
\draw (-1cm, -5.3cm) node[anchor=south,rotate=90] {BLIP+FINE};
\draw (-1cm, -7.45cm) node[anchor=south,rotate=90] {BLIP+C2F};

\end{tikzpicture}
\caption{Comparison of FINE with C2F using a constant initial guess (rows 2 and 3) and a good initial guess (rows 4 and 5) generated by BLIP with a coarse dictionary (row 1).}
\label{fig:brains}
\end{center}
\end{figure*}

\section{Conclusion}
\label{sec:conclusion}

In this paper, we have presented a temporal multiscale technique that approximates Bloch responses to reduce running times in qMRI. This was integrated into an optimisation framework that enabled the computation of cheap approximate gradients to be used in a reconstruction algorithm, which was found to improve the accuracy of solutions compared to similar methods without multiscaling. This ultimately means the same degree of accuracy can be achieved quicker. In contrast to dictionary-based methods e.g. BLIP, an optimisation approach does not suffer from memory problems or grid bias. Instead, an arbitrary degree of accuracy can be achieved. 

Areas for further research include providing error bounds for the multiscale Bloch responses, i.e. depending on the grid, by how much does the approximation differ to the true Bloch mapping. Another direction is to investigate the effect of varying the grid in the C2F algorithm. We only considered uniform grids but this is by no means necessary. Finally, it would be useful to find good strategies for choosing the vectors $N$ and $K$ in C2F and also criteria for when to refine the grid so that this can be done automatically.

\section{Compliance with Ethical Standards}
\label{sec:ethics}

This is a numerical simulation study for which no ethical approval was required.

\section{Acknowledgements}
\label{sec:acknowledgements}

S. Cortinhas acknowledges support from the Institute for Mathematical Innovation and the London Mathematical Society for undergraduate research bursaries. M. J. Ehrhardt acknowledges support from the EPSRC (EP/S026045/1, EP/T026693/1), the Faraday Institution (EP/T007745/1) and the Leverhulme Trust (ECF-2019-478).

\bibliographystyle{IEEEbib}
\bibliography{root}

\begin{thebibliography}{1}

\bibitem{bloch1946nuclear}
F.~Bloch,
\newblock ``Nuclear induction,''
\newblock {\em Physical Review}, vol. 70, no. 7-8, pp. 460--474, 1946.

\bibitem{ma2013magnetic}
D.~Ma, V.~Gulani, N.~Seiberlich, K.~Liu, J.~L. Sunshine, J.~L. Duerk, and M.~A.
  Griswold,
\newblock ``Magnetic resonance fingerprinting,''
\newblock {\em Nature}, vol. 495, no. 7440, pp. 187--192, 2013.

\bibitem{davies2014compressed}
M.~Davies, G.~Puy, P.~Vandergheynst, and Y.~Wiaux,
\newblock ``A compressed sensing framework for magnetic resonance
  fingerprinting,''
\newblock {\em SIAM Journal on Imaging Sciences}, vol. 7, no. 4, pp.
  2623--2656, 2014.

\bibitem{mcgivney2014svd}
D.~F. McGivney, E.~Pierre, D.~Ma, Y.~Jiang, H.~Saybasili, V.~Gulani, and M.~A.
  Griswold,
\newblock ``Svd compression for magnetic resonance fingerprinting in the time
  domain,''
\newblock {\em IEEE Transactions on Medical Imaging}, vol. 33, no. 12, pp.
  2311--2322, 2014.

\bibitem{asslander2018low}
J.~Assl{\"a}nder, M.~A. Cloos, F.~Knoll, D.~K. Sodickson, J.~Hennig, and
  R.~Lattanzi,
\newblock ``Low rank alternating direction method of multipliers reconstruction
  for mr fingerprinting,''
\newblock {\em Magnetic Resonance in Medicine}, vol. 79, no. 1, pp. 83--96,
  2018.

\bibitem{dong2019quantitative}
G.~Dong, M.~Hinterm{\"u}ller, and K.~Papafitsoros,
\newblock ``Quantitative magnetic resonance imaging: From fingerprinting to
  integrated physics-based models,''
\newblock {\em SIAM Journal on Imaging Sciences}, vol. 12, no. 2, pp. 927--971,
  2019.

\bibitem{sbrizzi2018fast}
A.~Sbrizzi, O.~van~der Heide, M.~Cloos, A.~van~der Toorn, H.~Hoogduin, P.~R.
  Luijten, and C.~A.~T. van~den Berg,
\newblock ``Fast quantitative mri as a nonlinear tomography problem,''
\newblock {\em Magnetic Resonance Imaging}, vol. 46, pp. 56--63, 2018.

\bibitem{scheffler1999pictorial}
K.~Scheffler,
\newblock ``A pictorial description of steady-states in rapid magnetic
  resonance imaging,''
\newblock {\em Concepts in Magnetic Resonance: An Educational Journal}, vol.
  11, no. 5, pp. 291--304, 1999.

\bibitem{golbabaee2017cover}
M.~Golbabaee, Z.~Chen, Y.~Wiaux, and M.~E. Davies,
\newblock ``Cover tree compressed sensing for fast mr fingerprint recovery,''
\newblock in {\em 2017 IEEE 27th International Workshop on Machine Learning for
  Signal Processing (MLSP)}. IEEE, 2017, pp. 1--6.

\end{thebibliography}

\end{document}